\documentclass[12pt]{article}

\usepackage[utf8]{inputenc} % allow utf-8 input
\usepackage{hyperref}       % hyperlinks
\usepackage{url}            % simple URL typesetting
\usepackage{booktabs}       % professional-quality tables
\usepackage{amsfonts}       % blackboard math symbols
\usepackage{nicefrac}       % compact symbols for 1/2, etc.
\usepackage{microtype}      % microtypography
\usepackage{xcolor}         % colors
\usepackage{algorithm}
\usepackage{algpseudocode}
\usepackage{subcaption}
\usepackage{amsmath}
\usepackage{comment}
\usepackage{natbib}
\usepackage{enumitem}
\usepackage{graphicx}
\usepackage{microtype}
\usepackage{multirow}
\usepackage{amsmath,amssymb,amsthm,mathtools,bm}
\usepackage{graphicx}
\usepackage{cleveref}
\usepackage{booktabs}     % colors
\newtheorem{theorem}{Theorem}

\newtheorem{definition}{Definition}
\newtheorem{assumption}{Assumption}

\newtheorem{remark}{Remark}

\def\o{{\scriptstyle{\mathcal{O}}}}
\def\O{\mathcal{O}}
\def\P{\mathbb{P}}

\newcommand{\E}{\mathbb{E}}
\newcommand{\PP}{\mathbb{P}}
\newcommand{\Var}{\mathrm{Var}}

\newcommand{\1}{\mathbf{1}}

\begin{document}

\def\spacingset#1{\renewcommand{\baselinestretch}%
{#1}\small\normalsize} \spacingset{1}

\date{}
%\if0\blind
%{
\title{\bf Asymptotic Anytime-Valid Inference for U-statistics}
\author{Leheng Cai\footnotemark[1] \footnotemark[3], $\,\,$  Qirui Hu\footnotemark[2] \footnotemark[3] \footnotemark[4]$\,\,$  and$\,\,$ Weijia Li\footnotemark[1]  \footnotemark[3] \hspace{.2cm}}
\maketitle
\renewcommand{\thefootnote}{\fnsymbol{footnote}}

\footnotetext[1]{Tsinghua University}
\footnotetext[2]{Shanghai University of Finance and Economics}
 \footnotetext[3]{ The authors are listed in alphabetical order with equal contribution}
\footnotetext[4]{Corresponding author: huqirui@mail.shufe.edu.cn}
%}

\begin{abstract}
We study asymptotic anytime-valid confidence sequences  for degree-two U-statistics under continuous monitoring. In the nondegenerate case, Hoeffding's projection reduces the problem to a time-uniform central limit theory for the partial sums of the  first-order projection, while the canonical remainder is shown to be negligible under mild moment assumptions. A leave-one-out jackknife  estimator then yields a fully data-driven procedure, leading to  confidence sequences  with asymptotic coverage guarantee for the parameter of interest. In the degenerate case, we show that the U-statistic is approximated by a centered quadratic Gaussian-chaos rather than by a simple Gaussian, which poses significant challenges for sequential inference. To address this issue, we novelly develop the Spectrally Allocated Gaussian-chaos Excursion (SAGE) boundary, and then provide plug-in implementations based on truncated  spectrum estimation with consistency guarantees. The resulting widths can attain the expected time-uniform optimal rates: \( \sqrt{\log\log n/n}\) in the nondegenerate regime and \( \log\log n/n\) in the degenerate regime. Several widely used U-statistics are discussed within the proposed framework, and  numerical experiments further support the validity of the derived theory.
\end{abstract}

\section{Introduction}
U-statistics are a typical device for estimating distributional functionals that can be written as expectations of symmetric kernels.  Since Hoeffding's decomposition, they have served as a unifying language for rank statistics, variance and covariance estimation, measures of dispersion and association, and two-sample or independence testing \citep{Hoeffding1948,serfling2009approximation,Lee1990,KorolyukBorovskich1994,MannWhitney1947,Kendall1938,sgaU2026}.  Their appeal is not only unbiasedness: the same projection calculus gives classical limit theory and separates the regular first-order component from genuinely second-order behavior.  In the nondegenerate regime, the first Hoeffding projection reduces the leading term to an average of i.i.d. influence functions, while in the canonical degenerate regime, the  projection vanishes and the statistic is  instead   a quadratic form or Gaussian chaos, with limit laws depending on the spectrum of the associated integral operator \citep{Hall1979,DehlingDenkerPhilipp1984,deJong1987,ArconesGine1993,deLaPenaGine1999,Gine2001lil,Ferger2003}.

This dichotomy is especially relevant in modern machine learning and nonparametric distributional inference.  Kernel maximum mean discrepancy (MMD), Hilbert--Schmidt independence criteria, energy distances, and distance covariance are all naturally expressed through U- or V-statistics; under the null hypothesis of equality or independence, their first-order projections typically vanish, producing the   weighted chi-square or Gaussian-chaos null limits \citep{GrettonBorgwardtRaschScholkopfSmola2012,BaringhausFranz2004,SzekelyRizzoBakirov2007,SzekelyRizzo2013,SejdinovicSriperumbudurGrettonFukumizu2013}.  Thus, the degenerate case is not a pathological corner case.  It is the operating regime of many distributional tests.  At the same time, it is the regime in which calibration is most delicate, because the boundary depends on spectral quantities that are usually unknown and must be estimated from the same data stream \citep{KoltchinskiiGine2000,RosascoBelkinDeVito2010}.

Anytime-valid inference asks for uncertainty statements that remain valid under continuous monitoring and optional stopping.  For sample means and related estimating equations, this goal is now well understood through nonnegative martingales, test supermartingales, mixture boundaries, and betting constructions \citep{Ville1939,DarlingRobbins1967,robbins1970statistical,robbins1970boundary,HowardRamdasMcAuliffeSekhon2021,betting,RamdasGrunwaldVovkShafer2023}.  These methods give finite-sample confidence sequences (CSs) when enough structure is available, for example boundedness, sub-exponential tails, or an exact martingale estimating equation.  A complementary line of work develops asymptotic confidence sequences (AsympCSs) by replacing fixed-time central limit approximations with time-uniform strong approximations, thereby recovering CLT-like versatility in sequential settings \citep{WaudbySmithArbourSinhaKennedyRamdas2024,GnettnerKirch2025}.

For U-statistics, however, the usual mean-based sequential tools do not transfer directly.  Adding a new observation changes all pairs involving that observation, so the statistic is not a sum of independent increments.  Even in the nondegenerate case, the reduction to a sample mean is only asymptotic and must be uniform over all monitoring times; a fixed-time CLT does not control the event of ever crossing a boundary.  One must prove that the canonical Hoeffding remainder is almost surely negligible relative to the iterated-logarithm order, and any plug-in  estimation does not disturb time-uniform coverage. Notably, the degenerate case is much harder: there is no Brownian first-order limit, the correct width is of order \(\log\log n/n\) rather than \(\sqrt{\log\log n/n}\), and a boundary must control a countable quadratic Gaussian chaos with unknown signed eigenvalues.

There is a substantial sequential literature around U-statistics, but it addresses different inferential targets.  Classical work studied fixed-width or two-stage sequential confidence procedures and stopping times based on U-statistics \citep{Sproule1985,MukhopadhyayVik1985,MukhopadhyayVik1988,Nasari2009}.  More recent work has used U-statistic processes for retrospective or online change-point detection, including robust Wilcoxon-type monitoring and distributional monitoring based on degenerate U-statistics \citep{DehlingFriedGarciaWendler2015,DehlingVukWendler2022,KirchStoehr2022,BonieceHorvathTrapani2025}.  These contributions are closely related in motivation, but their calibration is designed for testing change-points or for prescribed stopping rules.  They do not provide a general AsympCS for a fixed degree-two U-functional, and they do not give a unified treatment of the nondegenerate and canonical degenerate regimes with data-driven spectral calibration.

This paper fills that gap by developing asymptotic anytime-valid inference for degree-two U-statistics under continuous monitoring.  Our central idea is to separate the construction into two components: a probabilistic approximation that holds uniformly over time, and a statistical calibration step that estimates the nuisance quantities appearing in the limiting boundary.  This separation yields a modular framework and clarifies precisely where moment conditions, degeneracy, and spectral estimation enter.  In the nondegenerate regime, we show that the relevant limiting object is a Brownian partial-sum process associated with the first Hoeffding projection.  Based on an almost sure Gaussian approximation  and a strongly consistent jackknife variance estimator  (Theorems \ref{THM:GA1} and \ref{THM:strong_cons}), we construct AsympCSs; see Theorem \ref{THM:CS1}. These sequences connect naturally to classical time-uniform Gaussian boundaries; in particular, the stitched law-of-the-iterated-logarithm boundary achieves the optimal  width \( \sqrt{\log\log n/n}\).  
\par In the degenerate regime, the limiting object is instead a centered quadratic Gaussian-chaos process.  For this case, we derive a Gaussian-chaos approximation and novelly develop the Spectrally Allocated Gaussian-chaos Excursion (SAGE) boundary; see Theorems \ref{thm:degenerate_GA} and \ref{thm:gaussian_chaos}, which allocates the significance level budget across spectral coordinates and yields an AsympCS with the canonical \( \log\log n/n\) scale, matching the law-of-the-iterated-logarithm behavior of canonical degree-two U-statistics.  We then provide data-driven plug-in procedures for the degenerate regime based on truncated empirical spectral estimation. Finally, we illustrate the scope of the framework through   common examples, including sample variance, Gini mean difference, spatial Kendall's tau and  MMD. %and Appendix \ref{app:examples}. 
To our best knowledge, this paper first addressed the asymptotic anytime-valid inference for U-statistics systematically.  

The rest of the paper is organized as follows. Section~\ref{sec:preliminaries} fixes notations, recalls Hoeffding's decomposition, defines AsympCSs, and states the Gaussian boundaries used in the constructions.  Section~\ref{sec:theory} gives the nondegenerate and degenerate theoretical guarantees. Section~\ref{sec:examples} collects {\color{black}concrete examples of U-statistics.} Section~\ref{sec:simulations} contains simulation designs and numerical  results. Section~\ref{sec:conclusion} concludes. Detailed proofs and additional numerical results are deferred to the appendix. 

%{=========Prof.Cai==========}
\section{Preliminaries}
\label{sec:preliminaries}

%\subsection{Hoeffding decomposition}

\par \textbf{Hoeffding decomposition.} Let \(X_1,\ldots,X_n\) be i.i.d.\ random elements taking values in \(\mathcal X\). Consider the second-order \(U\)-statistic
\[
U_n
=
\binom{n}{2}^{-1}
\sum_{1\le i<j\le n}
h(X_i,X_j),
\]
where \(h:\mathcal X^2 \to \mathbb R\) is a symmetric kernel satisfying \(\mathbb E h^2(X_1,X_2)<\infty\). Let
$\theta
=
\mathbb E\,h(X_1,X_2)$ be the parameter of interest. Define the first-order projection
$
h_1(x)
=
\mathbb E\{h(x,X_2)\}-\theta$, and  the canonical kernel
$
h_2(x,y):=h(x,y)-\theta-h_1(x)-h_1(y).
$ 
%Then \(\E[h_1(X_1)]=0\), \(\E[h_2(X_1,X_2)]=0\), and \(\E[h_2(x,X_2)]=0\) for \(P\)-almost every \(x\). 
Hoeffding's decomposition gives
\begin{equation}
U_n-\theta
=
\frac{2}{n}\sum_{i=1}^n h_1(X_i)+R_n,
\qquad
R_n:=\binom{n}{2}^{-1}\sum_{1\le i<j\le n} h_2(X_i,X_j).
\label{eq:hoeffding-main}
\end{equation}
We call the U-statistic nondegenerate if \(\Var(h_1(X_1))>0\), and degenerate (canonical) if \(h_1(X_1)=0\) almost surely. For degenerate \(U\)-statistics, let \(\mathcal K:L^2(\P)\to L^2(\P)\) be the Hilbert-Schmidt operator
$(\mathcal K f)(x)=\int \{h(x,y)-\theta\}f(y) d\P(y)$. 
Then \(\mathcal K\) is compact and self-adjoint, and admits the spectral expansion in $L^2(\P\times \P)$, i.e. for  the eigensystems \((\lambda_\ell,\psi_\ell)_{\ell\ge1}\)  of \(\mathcal K\),
\begin{align}\label{eigenvalue}
    h (x,y) -\theta
=
\sum_{\ell\ge1}\lambda_\ell \psi_\ell(x)\psi_\ell(y). 
\end{align}

%{introduce the KL expansion...}

% Let \(X_1,X_2,\dots\) be i.i.d.\ random variables on a measurable space \((\mathcal X,\mathcal A)\), and let \(h:\mathcal X\times\mathcal X\to\R\) be a measurable symmetric kernel. The parameter of interest is
% $
% \theta := \E[h(X_1,X_2)],
% $
% and the degree-two U-statistic is
% \begin{equation}
% U_n := \binom{n}{2}^{-1}\sum_{1\le i<j\le n} h(X_i,X_j), \qquad n\ge 2.
% \end{equation}
% Write
% $
% m(x):=\E[h(x,X_2)], \qquad h_1(x):=m(x)-\theta,
% $
%\subsection{Asymptotic anytime-valid sequential inference}

% The Hoeffding decomposition in \eqref{eq:hoeffding-main} reduces the leading behavior of a nondegenerate second-order \(U\)-statistic to an average of i.i.d. first-order projections. For fixed-time inference, this representation naturally leads to classical central limit theory. Sequential inference, however, requires a stronger form of approximation: the resulting confidence set must remain valid uniformly over an infinite time horizon, allowing the stopping time to be chosen adaptively from the data.

 \noindent \textbf{Asymptotic anytime-valid sequential inference.} In this paper, we adopt the following framework of asymptotic anytime-valid  inference \citep{WaudbySmithArbourSinhaKennedyRamdas2024}. 
\begin{definition}[AsympCSs]
    Let $\mathcal{T}$ be a totally ordered infinite set. The intervals $\left[\widehat{\theta}_t-\ell_t, \widehat{\theta}_t+u_t\right]_{t \in \mathcal{T}}$ 
    %centered at the estimators $\left\{\widehat{\theta}_t\right\}_{t \in \mathcal{T}}$ 
    with $\ell_t, u_t>0$ form an $(1-\alpha)$-AsympCS for a sequence of real parameters $\left\{\theta_t\right\}_{t \in \mathcal{T}}$ if there exists a (typically unknown) non-asymptotic $(1-\alpha)$-CS  $\left[\widehat{\theta}_t-\ell_t^{\star}, \widehat{\theta}_t+u_t^{\star}\right]_{t \in \mathcal{T}}$ such that the probability
$%\begin{align*}
   \mathbb{P}\left(\forall t \in \mathcal{T}, \theta_t \in\left[\widehat{\theta}_t-\ell_t^{\star}, \widehat{\theta}_t+u_t^{\star}\right]\right) \geqslant 1-\alpha,
$%\end{align*}
~and $\ell_t^{\star} / \ell_t \xrightarrow{\text{a.s.}} 1$,\, $u_t^{\star} / u_t \xrightarrow{\text{a.s.}} 1$.
\end{definition}
Unlike existing nonasymptotic inference results for means of random variables, which often rely on strong boundedness or moment-generating-function assumptions \citep{betting,HowardRamdasMcAuliffeSekhon2021}, this asymptotic framework allows us to establish universal closed-form   boundaries for both nondegenerate and degenerate \(U\)-statistics under mild moment conditions.  

% The Hoeffding decomposition in \eqref{eq:hoeffding-main} reduces the leading behavior of a nondegenerate second-order \(U\)-statistic to an average of i.i.d. first-order projections. For classical fixed-time asymptotic inference, this representation naturally leads to confidence intervals whose critical values are determined by standard normal quantiles via the central limit theorem.

On the other hand,  a central tool for constructing AsympCSs  for the mean of i.i.d. random variables is a time-uniform boundary for Gaussian partial sums \citep{WaudbySmithArbourSinhaKennedyRamdas2024}, of which the formal definition is given below. %In our setting, once the \(U\)-statistic is approximated by a Gaussian process uniformly over time, a Gaussian boundary can be transferred to the statistic after suitable normalization. This separates the probabilistic part of the construction, which controls Gaussian fluctuations uniformly over time, from the statistical part, which establishes the required strong approximation and variance estimation.

\begin{definition}[Gaussian boundary]\label{def:boundary}
Fix $\alpha\in(0,1)$. A map $\gamma_{\alpha,m} :\mathbb N\to(0,\infty)$    is called a $(1-\alpha)$ Gaussian boundary starting from $m$ if, for i.i.d. standard Gaussian random variables $\{Z_i\}_{i\in\mathbb N}$,
\begin{equation}\label{eq:boundary_def}
  \P\left(\forall n\ge m:\ \left|\frac{1}{n}\sum_{i=1}^n Z_i\right|\le \gamma_{\alpha,m}(n)\right)\ge1-\alpha.
\end{equation}
\end{definition}
The starting point \(m\) marks the end of the cold-start phase.  %Different choices of \(\gamma_{\alpha,m}\) lead to different trade-offs between finite-sample width, asymptotic rate, and implementation simplicity. 
Several standard choices of the boundary \(\gamma_{\alpha,m}(\cdot)\) are available for non-asymptotic CSs of Gaussian means in the literature.  A representative example is the stitched boundary of \citep{HowardRamdasMcAuliffeSekhon2021}, which achieves the concentration rate \( \O\!\left(\sqrt{\log \log n / n}\right) \):
\begin{align}\label{LIL}
  \gamma_{\alpha,m}^{\text{LIL}}(n)
    =
  \frac{\eta^{1/4}+\eta^{-1/4}}{\sqrt{2n}} \sqrt{ {s\log\log\left(\max\left\{ \frac{\eta n}{m},e\right\}\right)+\log \frac{\zeta(s)}{\alpha (\log\eta)^s}  } },
\end{align}
where $\eta,s>1$, and $\zeta(\cdot)$ is the Riemann zeta function. 

Another classical construction is Robbins' Gaussian mixture boundary \citep{robbins1970statistical,robbins1970boundary}, which attains the rate \( \O(\sqrt{\log n / n}) \):
\begin{align}\label{GM}
    \gamma_{\alpha,m}^{\text{GM}}(n)
    =
    \sqrt{ \left[{\{g^{-1}(\alpha)\}^2+\log(n/m)}\right]/{n}},
    \quad \text{where} \;\;
    g(a)=2\{1-\Phi(a)+a\phi(a)\}.
\end{align}
Here, \(\Phi(\cdot)\) and \(\phi(\cdot)\) denote the CDF and PDF of a standard Gaussian distribution, respectively.
These bounds will be used as building blocks for our proposed anytime-valid inference procedure for U-statistics.

\section{Main results}
\label{sec:theory}
\par {\textbf{{Non-degenerate case.}}} We first consider that the kernel is non-degenerate at the first order, i.e.,
$\sigma^2:=\Var\{h_1(X_1)\}>0$. 

The Hoeffding decomposition in \eqref{eq:hoeffding-main} implies that the leading behavior of a nondegenerate second-order \(U\)-statistic is governed by an average of i.i.d. first-order projections. 
By establishing an almost sure Gaussian approximation for the first Hoeffding projection with positive variance and showing that the canonical remainder is asymptotically negligible almost surely, Theorem~\ref{THM:GA1} rigorously shows that the nondegenerate \(U\)-statistic can be approximated almost surely by the mean of a sequence of i.i.d. Gaussian random variables.
\begin{theorem}  \label{THM:GA1}
Suppose that \(\mathbb E|h(X_1,X_2)|^{2+\delta}<\infty\) for some \(\delta>0\), and that the kernel is non-degenerate at the first order, i.e., \(\sigma^2>0\).
On an enriched probability space there exists i.i.d.\ standard Gaussian random variables \(Z_1,\ldots,Z_n\)  such that \begin{align*}
    \left|U_n-\theta-\frac{2\sigma}{n}\sum_{i=1}^n Z_i\right|
=
\o_{a.s.}\!\left(n^{-1+1/(2+\delta)}\right).
\end{align*}
% \[
% \left|U_n-\theta-\frac{2\sigma}{n}\mathcal W(n)\right|
% =
% \o_{a.s.}\!\left(n^{-1/(2+\delta)}\right).
% \]
% Moreover, the leading term admits the distributional representation
% \[
% \frac{2\sigma}{n}\mathcal W(n)
% \;\overset{d}{=}\;
% \frac{2\sigma}{n}\sum_{i=1}^n Z_i.
% \]
% where \(Z_1,\ldots,Z_n\) are i.i.d.\ standard Gaussian random variables.
\end{theorem}
%The proof consists of two parts: establishing an almost sure Gaussian approximation for the first Hoeffding projection with positive variance, and showing that the canonical remainder is asymptotically negligible almost surely.

% Define \(\widehat h_{1,n}(X_i)=A_{i,n}-U_n\), where
% \[
% A_{i,n}
% :=
% \frac{1}{n-1}\sum_{j\ne i} h(X_i,X_j),
% \qquad
% 1\le i\le n.
% \]

We consider the jackknife-type variance estimator for the unknown variance \(\sigma^2\)  in Theorem~\ref{THM:GA1}, 
\begin{align}\label{sigma_hat}
  \widehat\sigma_n^2 =
\frac1n\sum_{i=1}^n \left\{\frac{1}{n-1}\sum_{j\ne i} h(X_i,X_j)\right\}^2-U_n^2.
\end{align}
The following Theorem \ref{THM:strong_cons} establishes the strong consistency of $\widehat\sigma_n^2$. 
\begin{theorem}\label{THM:strong_cons}
Suppose that \(\mathbb E|h(X_1,X_2)|^{4}<\infty\). The variance estimator $\widehat\sigma^2_n$ defined in (\ref{sigma_hat}) is strongly consistent:  
     $
\left|\widehat\sigma_n^2
-\sigma^2\right|
=
 \O_{a.s.}\left(\sqrt{ {\log\log n}/{n}}\right).
$
\end{theorem}
%Therefore, \(\sigma^2\) can be estimated at a polynomial rate, rather than merely at an arbitrary \(\o_{a.s.}(1)\) rate. This quantitative strong consistency is crucial for establishing the asymptotic time-uniform coverage guarantee  (\ref{time-uniform_cover}) in Theorem~\ref{THM:CS1}.

% For \(m=1\), \citep{WaudbySmithArbourSinhaKennedyRamdas2024} extended the Gaussian mixture boundary to the following form:
% \begin{align}\label{EQ:CS**}
%     \gamma_{\alpha,1}(n)
%     =
%     \sqrt{\frac{2(n\rho^2+1)}{n^2\rho^2}
%     \log\!\left( \frac{\sqrt{n\rho^2+1}}{\alpha}\right)},
%     \quad \forall \rho>0.
% \end{align}
% Here, the hyperparameter \(\rho\) can be tuned to minimize the interval width at a given time point for a fixed significance level \(\alpha\). 

For classical   asymptotic inference for non-degenerate U-statistics, one typically applies the central limit theorem and uses non-asymptotic standard normal quantiles as critical values to construct asymptotic confidence intervals (CIs).  
Motivated by Theorem~\ref{THM:GA1} and the analogy with classical asymptotic inference, we use Gaussian boundaries  to construct AsympCSs. 
The following Theorem \ref{THM:CS1} provides the corresponding theoretical justification.
\begin{theorem}\label{THM:CS1}
    Suppose that $\E |h(X_1,X_2)|^4<\infty$. Then, $\left[
U_n \pm  2\widehat \sigma_n\gamma_{\alpha,m}(n)
\right]_{n\geq m}$ forms a $(1-\alpha)$-AsympCS for $ \theta$ starting from $m$. 
 Furthermore,  \begin{align}\label{time-uniform_cover}
    \liminf_{m\to\infty}\P\left\{\forall n\geq m: 
    \theta\in\left[
U_n \pm  2\widehat \sigma_n\gamma_{\alpha,m}(n)
\right]\right\} \geq 1-\alpha.
\end{align} 
\end{theorem}
% {plug-in two bounds. LIL achieves the optimal rate.}
By plugging in the Gaussian boundaries in \eqref{LIL} or \eqref{GM}, we obtain universal AsympCSs for nondegenerate \(U\)-statistics. In particular, the plug-in boundary based on \eqref{LIL} achieves the optimal rate, since the law of the iterated logarithm for nondegenerate \(U\)-statistics yields
\[
    \limsup_{n\to\infty}
    \sqrt{ {n}/{(2\log\log n)}}
    |U_n-\theta|
    =
    2\sigma,
    \quad \text{a.s.}
\] In addition, we note that as long as \(\sigma^2\) admits a polynomial-rate estimator rather than merely an arbitrary \(\o_{a.s.}(1)\)-consistent estimator, the resulting AsympCS enjoys the asymptotic time-uniform coverage guarantee in \eqref{time-uniform_cover}. 
%We note that the asymptotic time-uniform coverage property in \eqref{time-uniform_cover} holds as long as \(\sigma^2\) can be estimated at a polynomial rate, rather than merely at an arbitrary \(\o_{a.s.}(1)\) rate.
%\subsection{Degenerate regime}

\noindent {  \textbf{Degenerate case. 
}} We now consider the canonical case, where \(h_1(X_1)=0\) almost surely. 
In this regime, the \(U\)-statistic is no longer dominated by the first-order projection in the Hoeffding decomposition; instead, it is entirely   the second-order degenerate component. 
Consequently, its asymptotic behavior changes significantly. 

The following Theorem \ref{thm:degenerate_GA} shows that the degenerate \(U\)-statistic admits an almost sure approximation by a centered quadratic Gaussian chaos process.
%the Gaussian-chaos almost sure invariance principle. 

\begin{theorem}\label{thm:degenerate_GA}
Assume that $\mathbb{E}\left\{h(X,X')\log |h(X,X')|\right\}^2 + \mathbb{E}\left\{h(X,X)\log |h(X,X)|\right\}^2<\infty$ and the kernel is degenerate at the first order, i.e., $\sigma^2=0$. 
    On an enriched probability space there exist independent standard Brownian motions \(\{ W_\ell(\cdot)\}_{\ell\ge 1}\)  such that
one has
\begin{equation}
\left|
U_n-\theta -n^{-2}\sum_{\ell\ge 1}\lambda_\ell\left\{W^2_\ell(n)-n\right\}
\right|
=
\o_{a.s.}\left(\log\log n/n\right),
\label{eq:deg-asip-main}
\end{equation}
{\color{black}where $\{\lambda_\ell\}_{\ell\geq 1}$ are the eigenvalues  defined in \eqref{eigenvalue}.}
\end{theorem}

Let \((\beta_{\ell,n})_{\ell\ge 1}\) be a nonnegative weight array with   \(\sum_{\ell\ge 1}\beta_{\ell,n}= 1\), and  we write \(\beta_{\ell}\) for \(\beta_{\ell,n}\) for simplicity. Let $\Lambda = \sum_{\ell\geq1 } \lambda_\ell $, $\Lambda^+ = \sum_{\ell:\lambda_\ell>0 } \lambda_\ell$, $\Lambda_\beta^+ = \sum_{\ell:\lambda_\ell>0 } \lambda_\ell \log(1/\beta_{\ell})$ and $\Lambda_{\beta,g}^+=\sum_{\ell:\lambda_\ell>0} \lambda_\ell  \left\{g^{-1}(\alpha\beta_\ell)\right\}^2$,  where \(g(a)=2\{1-\Phi(a)+a\phi(a)\}\). We now introduce the Spectrally Allocated Gaussian-chaos Excursion (SAGE) boundary $\Upsilon_{\alpha,m}(\cdot)$, for controlling quadratic Gaussian-chaos processes. For a fixed significance level \(\alpha\in(0,1)\), we distribute the  budget $\alpha$ across the Brownian coordinates by setting \(\alpha_\ell=\alpha\beta_\ell\), apply a Gaussian mean boundary to each \(W_\ell(\cdot)\), combine the coordinatewise guarantees via a union bound, and aggregate the resulting bounds with the spectral weights \( \lambda_\ell \). 
This yields the following  time-uniform upper boundary for the centered quadratic Gaussian-chaos process.

\begin{theorem}\label{thm:gaussian_chaos}
 Suppose that $\sum_{\ell\geq 1}|\lambda_\ell|<\infty$, $\Lambda_\beta^+<\infty$  and \(\{ W_\ell(\cdot)\}_{\ell\ge 1}\) are independent standard Brownian motions.  
For any \(\alpha\in(0,1)\), \(m\ge 1\), 
\[
\mathbb{P}\left(
\forall n\ge m:\ 
n^{-2}\sum_{\ell\ge 1}\lambda_\ell\bigl\{W_\ell^2(n)-n\bigr\}
\le
\Upsilon_{\alpha,m}(n)
\right)
\ge 1-\alpha,
\]
where \(\Upsilon_{\alpha,m}(n)\) can be chosen in either of the following two ways:
    \begin{align}\label{LIL2}
        &\Upsilon_{\alpha,m}^{\text{LIL}}(n)=\frac{\left(\eta^{1/4}+\eta^{-1/4}\right)^2}{ 2n }\left[  \left\{ {s\log\log\left(\max\left\{ \frac{\eta n}{m},e\right\}\right)+\log \frac{\zeta(s)}{\alpha (\log\eta)^s}  } \right\}\Lambda^++ \Lambda_\beta^+ \right]-\frac{\Lambda}{n},
    \end{align}  
\begin{align}\label{GM2}
    &\Upsilon_{\alpha,m}^{\text{GM}}(n)=n^{-1}\left[  \Lambda^+\log\left(n/m\right)+\Lambda_{\beta,g}^+-\Lambda\right].
\end{align}
\end{theorem}

When applying the Gaussian  boundaries in \eqref{LIL} and \eqref{GM} coordinatewise to each \(W_\ell(\cdot)\), we obtain two variants of the SAGE boundary, \eqref{LIL2} and \eqref{GM2}, with widths of order \(\O(\log\log n/n)\) and \(\O(\log n/n)\), respectively. The rate of $\Upsilon_{\alpha,m}(\cdot)$ in (\ref{LIL2}) matches the law of the iterated logarithm rate \cite{Gine2001lil}, since       $\limsup_{n\to\infty}
     { {n}}
    |U_n-\theta|/{ \log\log n}
    <\infty$ almost surely for canonical \(U\)-statistics.
 
%{Comment on the optimal rate: the LIL for canonical U-statistics...}

Moreover, we note  that the SAGE boundary proposed in Theorem~\ref{thm:gaussian_chaos} is a one-sided upper boundary.
This is sufficient for many classical statistical inference based on degenerate \(U\)-statistics, including distributional tests based on MMD or energy distance, where evidence against the null is typically accumulated in the upper tail. 
We defer the lower boundary results to Section~\ref{sec:lower} and focus on the upper boundary throughout the paper and the subsequent example.

In practice, implementation of the proposed SAGE boundaries requires consistent eigenvalue estimation, which we formalize in Assumption~\ref{ass1}.
\begin{assumption}\label{ass1}
  For a small $\iota>0$,  there exist consistent estimators $\widehat\Lambda$, $\widehat\Lambda^+$, $\widehat\Lambda_\beta^+$ and $\widehat\Lambda_{\beta,g}^+$ such that $|\widehat\Lambda-\Lambda
|+|\widehat\Lambda^+-\Lambda^+| +|\widehat\Lambda_\beta^+-\Lambda_\beta^+| + |\widehat\Lambda_{\beta,g}^+-\Lambda_{\beta,g}^+| =  \o_{a.s.}(n^{-\iota}).$ %\begin{align*}
%        &\widehat\Lambda-\Lambda = \o_{a.s.}(n^{-\iota}),\quad \widehat\Lambda^+-\Lambda^+=\o_{a.s.}(n^{-\iota}),\\&\widehat\Lambda_\beta^+-\Lambda_\beta^+=\o_{a.s.}(n^{-\iota}),\quad  \widehat\Lambda_{\beta,g}^+-\Lambda_{\beta,g}^+=\o_{a.s.}(n^{-\iota})
%    \end{align*}
\end{assumption}

\begin{remark}\label{remark1}
   Assumption~\ref{ass1} is a high-level condition, and we provide an empirical eigenvalue procedure that satisfies this condition. 
Consider the centered Gram matrix
\begin{equation*}
\widehat K_n(i,j)
:=
h(X_{i},X_{j})-U_n,
\qquad 1\le i,j\le  n.
%\label{eq:split-gram}
\end{equation*}
One can let \((\widehat\lambda_{\ell,n})_{\ell\ge 1}\) be the eigenvalues of \(n^{-1}\widehat {\mathbf K}_n\), in which the matrix $\widehat {\mathbf K}_n=\{\widehat {  K}_n(i,j)\}_{i,j=1}^{n}$. For simplicity, we write \(\widehat\lambda_\ell\) for \(\widehat\lambda_{\ell,n}\) when no confusion arises. 
Define $\widehat\Lambda =n^{-1}\sum_{i=1}^{n}h(X_i,X_i)-U_n$, $\widehat \Lambda^+=\sum_{\ell:\widehat\lambda_\ell>0}^{L_n} \widehat\lambda_\ell$, $\widehat \Lambda_\beta^+=\sum_{\ell:\widehat\lambda_\ell>0}^{L_n} \widehat\lambda_\ell \log(1/\beta_\ell)$, and $\widehat\Lambda_{\beta,g}^+=\sum_{\ell:\widehat\lambda_\ell>0}^{L_n} \widehat\lambda_\ell  \left\{g^{-1}(\alpha\beta_\ell)\right\}^2$, where
the truncation number  $L_n\asymp n^{a}$ and $\beta_\ell\asymp \ell^{-b}$  for some $0<a<1/2$ and $b>1$. Results in  related empirical operator analyses  including \citep{KoltchinskiiGine2000,RosascoBelkinDeVito2010} and   Weyl inequality ensure  that Assumption \ref{ass1} holds under regularity conditions. Further details are provided in Section  \ref{sec:details}.
\end{remark}

\begin{theorem}\label{thm:sequential_test}
 Under Assumption \ref{ass1} and the conditions in Theorem \ref{thm:degenerate_GA}-\ref{thm:gaussian_chaos}, 
for any \(\alpha\in(0,1)\), \(m\ge 1\), 
 % $\left(-\infty,U_n+\widehat\Upsilon_{\alpha,m}(n)\right]_{n\geq m}$ is an $(1-\alpha)$ AsympCS of $\theta$ starting from $m$, 
 \[
\liminf_{m\to\infty}\mathbb{P}\left(
\forall n\ge m:\ 
U_n
\le \theta+
\widehat\Upsilon_{\alpha,m}(n)
\right)
\ge 1-\alpha,
\]
where \(\widehat\Upsilon_{\alpha,m}(n)\) can be chosen in either of the following two ways: 
    \begin{align}\label{eq:LIL3}
        &\widehat\Upsilon_{\alpha,m}^{\text{LIL}}(n)=\frac{\left(\eta^{1/4}+\eta^{-1/4}\right)^2}{ 2n }\left[  \left\{ {s\log\log\left(\max\left\{ \frac{\eta n}{m},e\right\}\right)+\log \frac{\zeta(s)}{\alpha (\log\eta)^s}  } \right\}\widehat \Lambda^+ +\widehat \Lambda_\beta^+  \right]-\frac{\widehat\Lambda}{n},
    \end{align}
    \begin{align}\label{eq:GM3}
    &\widehat\Upsilon_{\alpha,m}^{\text{GM}}(n)=\frac{1}{n}\left[  \widehat\Lambda^+ \log\left(\frac{n}{m}\right)+\widehat\Lambda_{\beta,g}^+ -\widehat\Lambda\right].
\end{align}  
\end{theorem}

{\color{black}  Theorem \ref{thm:sequential_test} provides the asymptotic coverage guarantee for the plug-in SAGE boundaries.  The duality between CSs and sequential tests makes our framework well adapted for modern sequential testing and continuous monitoring. Formally, an asymptotic $(1 - \alpha)$-CS $\{\mathcal C_n\}_{n \ge m}$ is equivalent to a family of asymptotic level-$\alpha$ sequential tests for $H_0 : \theta = \theta_0$. By defining the rejection rule as $\phi_n = \mathbf{1}\{\theta_0 \notin \mathcal C_n\}$, the  Type I error is  asymptotically controlled at level $\alpha$ for any stopping time.  
%For instance, in a distribution-free A/B test based on the Kernel MMD, the null hypothesis $H_0:P=Q$ corresponds to the degenerate case where the target parameter $\theta=0$. By monitoring whether the SAGE-based confidence sequence excludes zero, we can terminate the trial early as soon as a treatment effect is detected.
}

{\par \noindent \textbf{Implementation. 
}} Building upon theoretical guarantees established for both regimes, we provide two algorithmic templates for constructing these anytime-valid AsympCSs below:

\begin{algorithm}[h]
\caption{Nondegenerate AsympCSs}
\label{alg:non-degenerate}
\begin{algorithmic}[1]
\Require Data stream \(X_1,X_2,\ldots\); symmetric kernel \(h\); level \(\alpha\in(0,1)\); delayed start \(m\ge2\). 
\Ensure Two-sided  \((1-\alpha)\)-AsympCSs \((\mathcal C_n^{\mathrm{ND}})_{n\ge m}\) for \(\theta=\E h(X_1,X_2)\).
\For{\(n=m,m+1,\ldots\)}
    \State 
    Compute
    \(U_n\gets \binom{n}{2}^{-1}\sum_{1\le i<j\le n}h(X_i,X_j)\).
    \State 
    Compute boundary \(\gamma_{\alpha,m}(\cdot)\) using (\ref{LIL}) or (\ref{GM}) and estimate \(\widehat\sigma_n^2\) from \eqref{sigma_hat}.

    \State 
    Output
    \(\mathcal C_n^{\mathrm{ND}}\gets [U_n-2\widehat\sigma_n\gamma_{\alpha,m}(n),\,U_n+2\widehat\sigma_n\gamma_{\alpha,m}(n)]\).
\EndFor
\end{algorithmic}
\end{algorithm}

\begin{algorithm}[h]
\caption{Degenerate  AsympCSs}
\label{alg:degenerate}
\begin{algorithmic}[1]
\Require Data stream \(X_1,X_2,\ldots\); canonical symmetric kernel \(h\); level \(\alpha\in(0,1)\); delayed start \(m\ge2\); truncation number \(L_n\);   weights \((\beta_\ell)_{\ell\ge1}\) with \(\sum_{\ell\ge1}\beta_\ell=1\)
\Ensure One-sided   \((1-\alpha)\)-AsympCSs \((\mathcal C_n^{\mathrm{D}})_{n\ge m}\) for \(\theta=\E h(X_1,X_2)\).
\For{\(n=m,m+1,\ldots\)}
    \State %\textbf{Update statistic:} 
    Compute
    \(U_n\gets \binom{n}{2}^{-1}\sum_{1\le i<j\le n}h(X_i,X_j)\).
    \State %\textbf{Estimate nuisance:}
    Form the centered Gram matrix
    \(\widehat {\mathbf K}_n\gets h(X_i,X_j)-U_n\), \(1\le i,j\le n\). \State Retain the top $L_n$ eigenvalues of the scaled matrix \(n^{-1}\widehat {\mathbf K}_n\) in absolute value: $\widehat{\lambda}_1, \dots, \widehat{\lambda}_{L_n}$.
    \State  %\textbf{Compute spectral penalty constants:}
    Compute $\widehat\Lambda$, $\widehat \Lambda^+$, $\widehat \Lambda_\beta^+$, $ \widehat\Lambda_{\beta,g}^+$ provided in Remark \ref{remark1}.
    \State %\textbf{Calibrate boundary:} 
    Compute the SAGE boundary \(\widehat\Upsilon_{\alpha,m}(n)\)   using \eqref{eq:LIL3} or \eqref{eq:GM3}.
    \State %\textbf{Report confidence set:}
    Output
    \(\mathcal C_n^{\mathrm{D}}\gets [U_n-\widehat\Upsilon_{\alpha,m}(n),\infty)\).
\EndFor
\end{algorithmic}
\end{algorithm}
%{Further details concerning the computational complexity of Algorithm \ref{alg:degenerate}, along with a scalable acceleration scheme, are discussed in Appendix \ref{app:comp-complexity}.} 
A key practical consideration in implementing Algorithm \ref{alg:degenerate} is the choice of the significance level allocation weights \((\beta_\ell)_{\ell\ge1}\). To tighten the  SAGE boundaries, a simple heuristic is to choose the allocation weights \((\beta_\ell)_{\ell\ge1}\) to minimize \(\Lambda_\beta^+\) under the constraints \(\beta_\ell\geq 0\) and \(\sum_{\ell\ge1}\beta_\ell=1\). 
Let $p_\ell= {\max\{\lambda_\ell,0\}}/{\Lambda^+}$. 
Then, 
\[
\Lambda_\beta^+
=
\sum_{\ell\geq 1}\max\{\lambda_\ell,0\}\log(1/\beta_\ell)
=
\Lambda^+ \left\{\sum_{\ell\ge1}p_\ell\log(1/p_\ell)+   \mathrm{KL}(\bm p  \|\bm \beta)\right\},
\]
where the KL divergence $\mathrm{KL}(\bm p\|\bm \beta)
:=
\sum_{\ell\ge1}p_\ell\log(p_\ell/\beta_\ell)$ is nonnegative and equals zero if and only if \(\beta_\ell=p_\ell\) for all \(\ell\ge1\). Hence, the optimal choice is
$\beta_\ell^\star
=
{\max\{\lambda_\ell,0\}}/{\Lambda^+}$, which motivates a data-driven choice  $\widehat\beta_\ell = \max\{\widehat\lambda,0\}/\widehat\Lambda^+$ where $\widehat\Lambda^+=\sum_{\ell:\widehat\lambda_\ell>0}\widehat\lambda_\ell$ for each $n=m,m+1,\ldots$.

%{To reduce the substantial $\mathcal{O}(n^3)$ per-step eigendecomposition cost in Algorithm 2, one can approximate the spectrum using an $N_n = \lceil n^\alpha \rceil$ subsample for $\alpha \in (0, 1)$. This scalable strategy reduces the computational complexity to $\mathcal{O}(n^{3\alpha})$ while rigorously preserving the asymptotic anytime-validity, provided that the subsampled approximation satisfies Assumption 1.}

% {\color{black}\begin{algorithm}[h]
% \caption{Data-driven weights $\{\beta_\ell\}_{\ell\geq 1}$ for the SAGE boundary}
% \label{alg:datadriven}
% \begin{algorithmic}[1]
% \Require  The top $L_n$ estimated eigenvalues in absolute value: $\widehat{\lambda}_1, \dots, \widehat{\lambda}_{L_n}$.
% \Ensure Plug-in optimal weights $\{\widehat\beta_\ell\}_{\ell\geq 1}$. 
% \For{\(n=m,m+1,\ldots\)}
%     \State Compute $\widehat\Lambda^+=\sum_{\ell:\widehat\lambda_\ell>0}\widehat\lambda_\ell$.
%     \State Output $\widehat\beta_\ell = \max\{\widehat\lambda,0\}/\widehat\Lambda^+$.
% \EndFor
% \end{algorithmic}
% \end{algorithm}}

%{Weijia: Discuss the data-driven approach. truncate... this part could be moved to appendix}

%Weijia: please discuss briefly how to turn the CS to sequential test.   put them in appendix if there is no space.}

%\clearpage
%{\color{red} Revised version. Please check them. ===============

%=============================
%}

%\clearpage
\section{Examples}
\label{sec:examples}

Here, we introduce four concrete examples, which are adopted in numerical experiments. 
\paragraph{Sample variance.}
The target parameter $\theta=\mathrm{Var}(X_1)$ is estimated by the kernel: $h(x,y) = (x-y)^2/2$. Let $\mu = \mathbb{E} X_1$. The first-order projection and its variance are given by:
\[
    h_1(x) =  \left\{(x-\mu)^2-\theta\right\}/2, \quad \sigma^2 = \mathrm{Var}\left\{(X_1-\mu)^2\right\}/4.
\]
If $X_1 \sim \mathcal{N}(\mu, \theta)$, the fourth central moment is $3\theta^2$, which yields a   closed-form variance $\sigma^2 = \theta^2/2$. 
Since  Theorem~\ref{THM:strong_cons} requires $\mathbb{E}h^4 < \infty$, the sample variance mandates an  eighth-moment condition on $X_1$, which explains its mild sensitivity to heavier-tailed data observed in our simulations. 

\paragraph{Gini mean difference (GMD).}
The target   $\theta = \mathbb{E}|X_1-X_2|$ is associated with the kernel: $h(x,y) = |x-y|$. The first-order projection and its corresponding variance are  defined as:
\[
    h_1(x) = \mathbb{E}|x-X_2| - \theta, \quad \sigma^2 = \mathrm{Var}\{\mathbb{E}(|X_1-X_2| \mid X_1)\}.
\]
If $X_1 \sim \mathcal{N}(\mu, v)$, the target evaluates to $\theta = 2\sqrt{v/\pi}$, and the projection variance admits the exact  form $\sigma^2
=
v\left\{1/3+(2\sqrt3-4)/\pi\right\}$. 
A finite  fourth moment  on $X_1$   guarantees the asymptotic validity, making GMD  more robust to heavy tails.

\paragraph{Spatial Kendall's tau.}
 Consider bivariate data $\bm X_i=(X_{i1},X_{i2})^\top  \in \mathbb{R}^2$. The kernel for the off-diagonal entry of the two-dimensional spatial Kendall's tau matrix  \citep{ChoiMarden1998} is:
\[
    h(\bm x,\bm y) = \frac{(x_1-y_1)(x_2-y_2)}{\|\bm x-\bm y\|_2^2} \mathbf{1}\{\bm x \neq \bm y\}.
\]
The target $\theta = \mathbb{E}\{h(\bm X_1, \bm X_2)\}\in[-1/2,1/2]$  provides a robust measure of directional association between two coordinates.   Under a bivariate elliptical distribution with shape correlation $\rho \in[-1,1]$, the target has a  closed form $\theta = \left(1-\sqrt{1-\rho^2}\right)/(2\rho)$ when $\rho\neq 0$. {\color{black}In our simulation, we set    $\rho=0.6$, which then gives $\theta=1/6$. {\color{black}Further details  are deferred to Appendix~\ref{app:addition-fig}.}}  Since the kernel is bounded, all moment conditions required by Theorem~\ref{THM:CS1} hold automatically, which explains the stable coverage observed in the heavy-tailed settings considered here.

\paragraph{MMD with Gaussian kernel.}
Let $\bm Z_i = (X_i, Y_i)^\top$ denote a paired observation where $X_i \sim P$ and $Y_i \sim Q$. The unbiased squared MMD statistic \citep{GrettonBorgwardtRaschScholkopfSmola2012} employs the kernel:
\[
    h(\bm z, \bm z') = k(x, x') + k(y, y') - k(x, y') - k(x', y),
\]
where $k(x,y) = \exp(-|x-y|^2/2)$ is the Gaussian kernel. The target parameter is $\theta = \E h(\bm Z_1,\bm Z_2)=:\mathrm{MMD}^2(P, Q)$. The first-order projection is defined generally as $h_1(\bm z) = \mathbb{E} h(\bm z, \bm Z)  - \theta$.  
Under the null hypothesis that $P=Q$, a  direct calculation yields $\theta = 0$ and  the projection vanishes:
$\sigma^2 = \mathrm{Var}\{h_1(\bm Z)\} = 0$, then the U-statistic completely loses its first-order term and falls exactly into the canonical degenerate case. 
 
\section{Simulations}
\label{sec:simulations}

\par \textbf{Non-degenerate case.} For non-degenerate U-statistics, we consider GMD under three data generating distributions: the standard Gaussian, $t_{10}$ and Laplace with density $f(x) =  \exp\left(-\sqrt{2}|x|\right)/\sqrt{2}$. To evaluate the   empirical coverage and the average width of  the proposed AsympCSs, we conduct Monte Carlo simulations with 500 independent replications.
We set the nominal error rate to \(\alpha = 0.05\) throughout, and monitor the U-statistics continuously from a   start $m=400$ up to $n_{\max} = 10,000$.

We compare against the baseline {\textbf{classical pointwise CIs}} based on the asymptotic normality and evaluate two variants of the AsympCS in
\eqref{time-uniform_cover}: {\textbf{AsympCS-LIL}} using the stitched boundary in 
\eqref{LIL} with \(\eta=2.0\) and \(s=1.4\) suggested by 
\cite{HowardRamdasMcAuliffeSekhon2021}, and {\textbf{AsympCS-GM}} using the
Gaussian-mixture boundary in \eqref{GM}. The results are illustrated in Figures \ref{fig:nondeg_main} and \ref{fig:la-gmd}-\ref{fig:t-gmd} in Appendix \ref{app:gmd}.
Additional simulation results for other U-statistics, including sample variance and  Spatial Kendall's tau, are deferred to Appendix~\ref{app:addition-fig}. As shown in the left panel of Figure~\ref{fig:nondeg_main}, the classical CIs suffer  from a severe loss of coverage due to repeated sequential evaluations; their cumulative miscoverage rate rapidly climbs    above the nominal $\alpha$ level. In contrast, both AsympCS-GM and AsympCS-LIL   maintain the miscoverage rate below $\alpha $ uniformly across the entire horizon. The plot of one sample path for $U_n$ in the right panel of Figure~\ref{fig:nondeg_main} further illustrates how the valid AsympCSs safely envelop the true parameter.

The middle panel of Figure~\ref{fig:nondeg_main} traces the average half-widths of the proposed AsympCSs for GMD under the standard Gaussian distribution. As expected from theory, the AsympCSs are wider than the classical pointwise CIs, reflecting the necessary correction for continuous monitoring. Notably, although the AsympCS-GM boundary shrinks at the slower rate \(\sqrt{\log n/n}\), it yields sharper intervals in the finite-sample situation. This observation is
consistent with \citep{HowardRamdasMcAuliffeSekhon2021}: the GM boundary can be tighter at moderate sample sizes, but is eventually outperformed by the
AsympCS-LIL boundary, which attains the asymptotically optimal \(\sqrt{\log\log n/n}\) rate. A similar phenomenon also emerges in our degenerate experiments, translating into higher finite-sample power using the SAGE-GM boundary. 

\begin{figure}[h] 
    \centering
    \includegraphics[width=\textwidth]{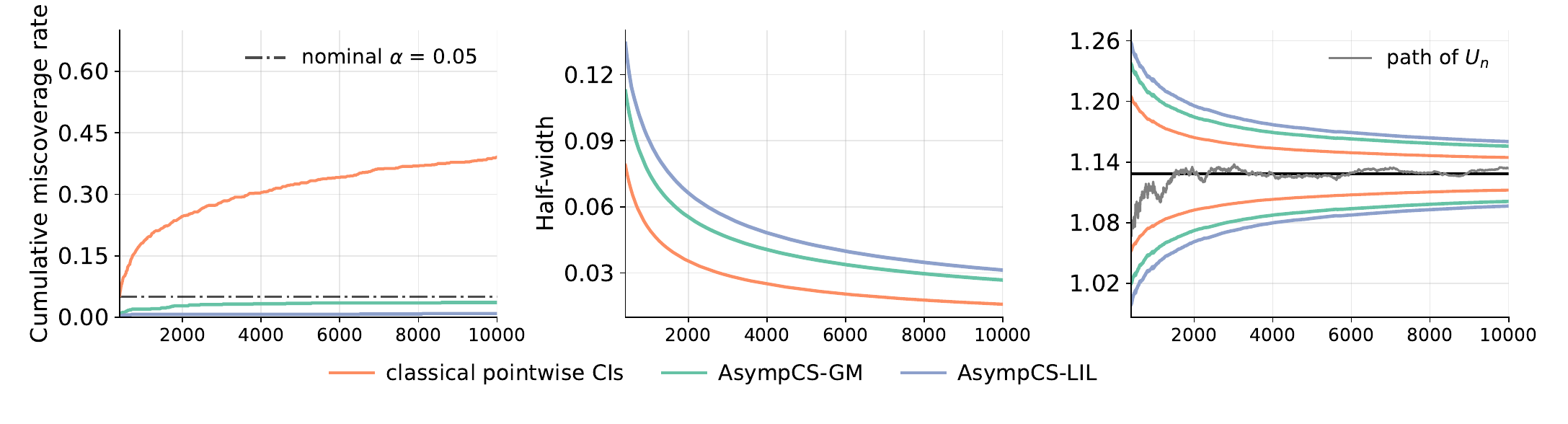} 
    \caption{{\color{black}Left: cumulative miscoverage rates of the proposed AsympCSs and the classical pointwise CIs for GMD under the standard Gaussian distribution. Middle: averaged half-widths of the proposed AsympCSs and the classical pointwise CIs. Right: a single sample path of the statistics $U_n$ for GMD alongside the three boundaries, where the black horizontal line indicates the true parameter. The horizontal axis in all three panels represents the sample size.}}
    \label{fig:nondeg_main}
\end{figure}
\begin{figure}[h]
    \centering
    \begin{subfigure}[b]{0.48\textwidth}
        \centering
        \includegraphics[width=\textwidth]{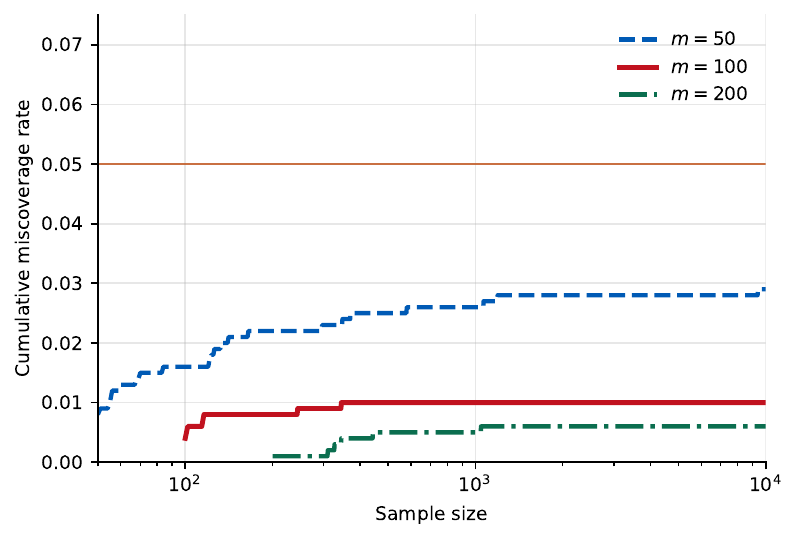}
        \label{fig:coldstart_gaussian}
    \end{subfigure}
    \hfill
    \begin{subfigure}[b]{0.48\textwidth}
        \centering
        \includegraphics[width=\textwidth]{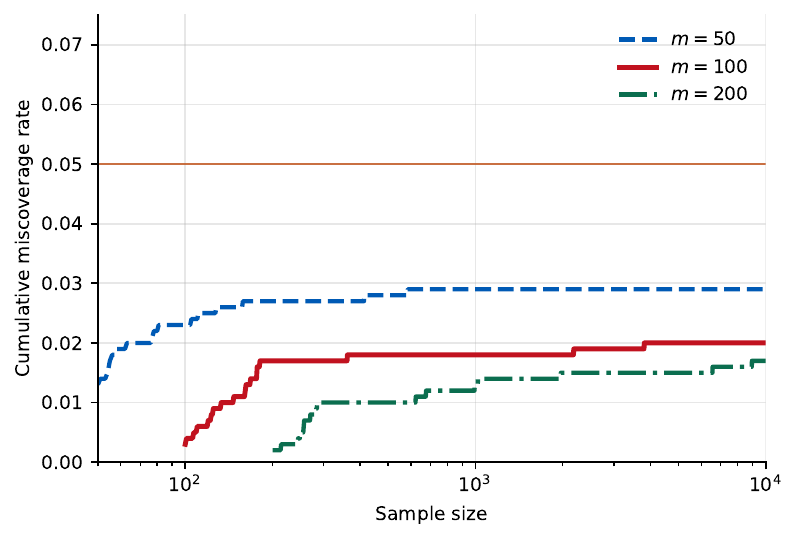}
        \label{fig:coldstart_t10}
    \end{subfigure}
    \caption{Sensitivity of AsympCS-LIL for GMD to the cold-start parameter $m\in\{50,100,200\}$. 
Cumulative miscoverage rates are reported under the standard Gaussian (left panel) and \(t_{10}\) (right panel) distributions.}
    \label{fig:coldstart_sensitivity}
\end{figure}

Also we consider the sensitivity to the cold-start parameter $m$. We examine the reliability of the asymptotic coverage under continuous monitoring by varying the cold-start parameter $m \in \{50, 100, 200\}$. 
Figure~\ref{fig:coldstart_sensitivity}  reports the cumulative miscoverage rates of the  AsympCS-LIL for  GMD under the standard Gaussian and the $t_{10}$ distributions. In both settings, once the start time   $m \ge 50$, the cumulative miscoverage rate remains below the level $\alpha$ for up to $10^4$ observations. 
It shows that a moderate cold-start absorbs early finite-sample deviations.

\par \noindent \textbf{Degenerate case.}
For degenerate U-statistics, we study   two-sample sequential testing   based on   MMD with   Gaussian kernel over the time horizon from \(m=400\) to \(n_{\max}=2{,}000\). { The truncation number $L_n$ is set to be $[n^{1/4}]$.} We evaluate the sequential MMD test under the standard Gaussian null \(H_0:P=Q\)  and under mean-shift alternatives \(H_1:P\neq Q\), parameterized by \(\delta\). We compare against the   baseline {\textbf{classical test}}   using  quantiles of a centered weighted sum of chi-squared random variables as critical values, and  employ the proposed SAGE boundaries in  \eqref{eq:LIL3} and  ~\eqref{eq:GM3}, denoted by {\textbf{SAGE-LIL}} and {\textbf{SAGE-GM}} respectively.  
All results are based on 500 independent Monte Carlo replications. Additional simulation results under other null distributions are deferred to Appendix~\ref{app:mmd}.

\begin{figure*}[h]
    \centering

    \includegraphics[width=\textwidth]{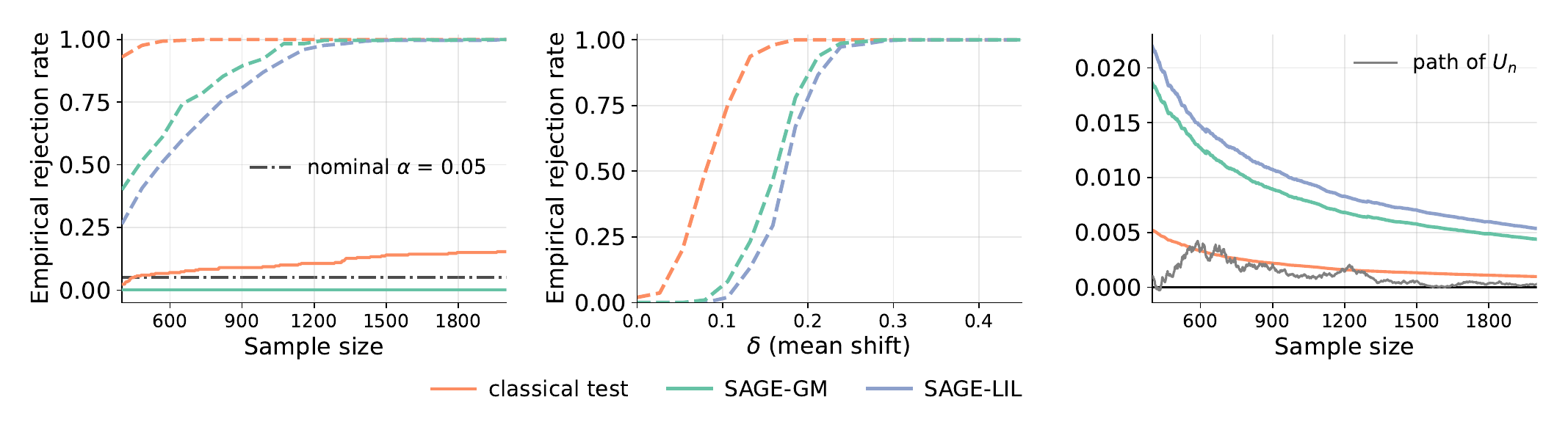}

    \caption{{\color{black}Left: the size and power comparison of the proposed sequential testing procedure using the SAGE boundaries and the classical test procedure for sequential two-sample test with MMD kernel statistics under standard Gaussian distribution. Dashed lines represent the power under $H_1$ with  $\delta = 0.3$, while solid lines represent the size under $H_0$ with $\delta = 0$. Middle: the empirical power   over mean shifts $\delta \in [0, 0.45]$ at sample size $n = 2000$. Right: a single sample path of the MMD test statistic $U_n$ alongside the three boundaries under $H_0$.}}
    \label{fig:deg_main}
\end{figure*}

As shown in the left panel of Figure~\ref{fig:deg_main}, the classical testing rule completely fails under continuous monitoring, with its empirical rejection rate rapidly accumulating beyond the nominal $\alpha$ level. 
In contrast, the SAGE boundaries yield sizes close to zero under \(H_0\), while their power approaches one as the sample size increases under \(H_1\) where $\delta=0.3$ in the simulation. It confirms that our framework achieves rigid Type I error control and has enough power simultaneously. 
The right panel of Figure~\ref{fig:deg_main} also visualizes the underlying mechanism via a single sample path under the null.

To examine the sensitivity to signal strength, we evaluate the power over a grid of mean shifts $\delta \in [0, 0.45]$ at the fixed sample size $n=2000$. As depicted in the middle panel of Figure~\ref{fig:deg_main}, the empirical power  based on our SAGE boundaries increases from zero to one, demonstrating their effectiveness for sequential testing.  

The empirical performance of the  sequential test   relies on the significance level budget  allocation across the infinite eigenvalue spectrum. We compared various allocation strategies, including polynomial decay {\color{black}{weights}} \(\beta_\ell=\ell^{-b}/\zeta(b)\), exponential decay {\color{black}{weights}} \(\beta_\ell=(1-e^{-c})e^{-c(\ell-1)}\), and our data-driven  {\color{black}{weights}} $\beta_\ell \propto \max\left\{\widehat{\lambda}_\ell,0\right\}/\widehat\Lambda^+$, as shown in Figure~\ref{fig:deg_sensitivity} in Appendix \ref{sec:sensitivity}. {\color{black}These results highlight the theoretical role of \(\log(1/\beta_\ell)\): 
a larger \(b\) concentrates the budget on the leading eigenvalues and thus widens 
the boundary, while a larger \(c\) inflates the boundary through the rapid tail 
growth of \(\log(1/\beta_\ell)\).} Remarkably, the proposed data-driven procedure optimally distributes the     budget, yielding the tightest  boundary while remaining below the nominal level. This justifies our adoption of this data-driven approach to compute $\beta_\ell$ throughout the degenerate sequential testing experiments.

\section{Concluding remarks}
\label{sec:conclusion}

This paper develops asymptotic anytime-valid inference for degree-two U-statistics under continuous monitoring. In the nondegenerate regime, Hoeffding's decomposition reduces the problem to the first-order projection, leading to AsympCS with LIL scale \(\sqrt{\log\log n/n}\). In the degenerate regime, the limit is instead a spectral quadratic Gaussian chaos with unknown signed eigenvalues. The proposed SAGE boundary separates probability calibration from spectral estimation and yields a one-sided plug-in AsympCSs with the sharper \(\log\log n/n\) scale under a weighted spectral approximation condition, covering degenerate kernel and distance-based tests.

Several limitations remain. The theory is restricted to scalar degree-two U-statistics; vector- or matrix-valued extensions require simultaneous time-uniform approximations and control of dimension-dependent errors, especially for spatial Kendall matrices and high-dimensional elliptical models \citep{ChoiMarden1998,DuerreTylerVogel2016,HanLiu2018}. The algorithms are not yet streaming-optimized, since exact updates, pairwise storage, and repeated eigendecompositions are costly; scalable versions may need incomplete U-statistics, incremental summaries, randomized low-rank updates, or linear-time kernel approximations \citep{GrettonBorgwardtRaschScholkopfSmola2012,ZarembaGrettonBlaschko2013,ChwialkowskiRamdasSejdinovicGretton2015}. Finally, the degenerate theory assumes a verified spectral approximation rate, and finite-sample cold starts, local degeneracy, and automatic switching between regimes remain open.

%\clearpage
{
\small
\bibliographystyle{plain}
\bibliography{ref}

@article{robbins1970boundary,
  author  = {Robbins, Herbert and Siegmund, David},
  title   = {Boundary Crossing Probabilities for the {Wiener} Process and Sample Sums},
  journal = {The Annals of Mathematical Statistics},
  year    = {1970},
  volume  = {41},
  number  = {5},
  pages   = {1410--1429},
  doi     = {10.1214/aoms/1177696787}
}

@article{robbins1970statistical,
  author  = {Robbins, Herbert},
  title   = {Statistical Methods Related to the Law of the Iterated Logarithm},
  journal = {The Annals of Mathematical Statistics},
  year    = {1970},
  volume  = {41},
  number  = {5},
  pages   = {1397--1409},
  doi     = {10.1214/aoms/1177696786}
}

@article{betting,
    author = {Waudby-Smith, Ian and Ramdas, Aaditya},
    title = {Estimating means of bounded random variables by betting},
    journal = {Journal of the Royal Statistical Society Series B: Statistical Methodology},
    volume = {86},
    number = {1},
    pages = {1-27},
    year = {2024},
    month = {02},
    issn = {1369-7412},
    doi = {10.1093/jrsssb/qkad009},
    url = {https://doi.org/10.1093/jrsssb/qkad009},
    eprint = {https://academic.oup.com/jrsssb/article-pdf/86/1/1/56961777/qkad009.pdf},
}

@book{serfling2009approximation,
  author    = {Serfling, Robert J.},
  title     = {Approximation Theorems of Mathematical Statistics},
  publisher = {John Wiley \& Sons},
  address   = {New York},
  year      = {1980},
  doi       = {10.1002/9780470316481}
}

@misc{BonieceHorvathTrapani2025,
  author       = {Boniece, Cooper and Horvath, Lajos and Trapani, Lorenzo},
  title        = {Sequential Monitoring for Distributional Changepoint Using Degenerate {U}-Statistics},
  year         = {2025},
  eprint       = {2510.22368},
  archivePrefix= {arXiv},
  primaryClass = {math.ST},
  doi          = {10.48550/arXiv.2510.22368}
}

@article{DarlingRobbins1967,
  author  = {Darling, D. A. and Robbins, Herbert},
  title   = {Iterated Logarithm Inequalities},
  journal = {Proceedings of the National Academy of Sciences of the United States of America},
  year    = {1967},
  volume  = {57},
  number  = {5},
  pages   = {1188--1192},
  doi     = {10.1073/pnas.57.5.1188}
}

@article{deJong1987,
  author  = {de Jong, Piet},
  title   = {A Central Limit Theorem for Generalized Quadratic Forms},
  journal = {Probability Theory and Related Fields},
  year    = {1987},
  volume  = {75},
  number  = {2},
  pages   = {261--277},
  doi     = {10.1007/BF00354037}
}

@book{deLaPenaGine1999,
  author    = {de~la Pe{\~n}a, Victor H. and Gin{\'e}, Evarist},
  title     = {Decoupling: From Dependence to Independence},
  publisher = {Springer},
  address   = {New York},
  year      = {1999},
  doi       = {10.1007/978-1-4612-0537-1}
}

@article{DehlingDenkerPhilipp1984,
  author  = {Dehling, Herold and Denker, Manfred and Philipp, Walter},
  title   = {Invariance Principles for von {Mises} and {U}-Statistics},
  journal = {Zeitschrift f{\"u}r Wahrscheinlichkeitstheorie und Verwandte Gebiete},
  year    = {1984},
  volume  = {67},
  number  = {2},
  pages   = {139--167},
  doi     = {10.1007/BF00535265}
}

@article{Ferger2003,
  author  = {Ferger, Dietmar},
  title   = {A Functional Law of the Iterated Logarithm for {U}-Statistic Type Processes},
  journal = {Acta Applicandae Mathematicae},
  year    = {2003},
  volume  = {78},
  pages   = {115--120},
  doi     = {10.1023/A:1025792223679}
}

@article{Gine2001lil,
  author  = {Gin{\'e}, Evarist and Kwapie{\'n}, Stanis{\l}aw and Lata{\l}a, Rafa{\l} and Zinn, Joel},
  title   = {The {LIL} for Canonical {U}-Statistics of Order 2},
  journal = {The Annals of Probability},
  year    = {2001},
  volume  = {29},
  number  = {1},
  pages   = {520--557},
  doi     = {10.1214/aop/1008956343}
}

@article{GnettnerKirch2025,
  author  = {Gnettner, Felix and Kirch, Claudia},
  title   = {A New and Flexible Class of Sharp Asymptotic Time-Uniform Confidence Sequences},
  journal = {Statistics \& Probability Letters},
  year    = {2025},
  volume  = {226},
  pages   = {110462},
  doi     = {10.1016/j.spl.2025.110462}
}

@article{GrettonBorgwardtRaschScholkopfSmola2012,
  author  = {Gretton, Arthur and Borgwardt, Karsten M. and Rasch, Malte J. and Sch{\"o}lkopf, Bernhard and Smola, Alexander J.},
  title   = {A Kernel Two-Sample Test},
  journal = {Journal of Machine Learning Research},
  year    = {2012},
  volume  = {13},
  pages   = {723--773}
}

@article{Hall1979,
  author  = {Hall, Peter},
  title   = {On the Invariance Principle for {U}-Statistics},
  journal = {Stochastic Processes and their Applications},
  year    = {1979},
  volume  = {9},
  number  = {2},
  pages   = {163--174},
  doi     = {10.1016/0304-4149(79)90015-2}
}

@article{Hoeffding1948,
  author  = {Hoeffding, Wassily},
  title   = {A Class of Statistics with Asymptotically Normal Distribution},
  journal = {The Annals of Mathematical Statistics},
  year    = {1948},
  volume  = {19},
  number  = {3},
  pages   = {293--325},
  doi     = {10.1214/aoms/1177730196}
}

@article{HowardRamdasMcAuliffeSekhon2021,
  author  = {Howard, Steven R. and Ramdas, Aaditya and McAuliffe, Jon and Sekhon, Jasjeet},
  title   = {Time-Uniform, Nonparametric, Nonasymptotic Confidence Sequences},
  journal = {The Annals of Statistics},
  year    = {2021},
  volume  = {49},
  number  = {2},
  pages   = {1055--1080},
  doi     = {10.1214/20-AOS1991}
}

@article{Kendall1938,
  author  = {Kendall, Maurice G.},
  title   = {A New Measure of Rank Correlation},
  journal = {Biometrika},
  year    = {1938},
  volume  = {30},
  number  = {1/2},
  pages   = {81--93},
  doi     = {10.1093/biomet/30.1-2.81}
}

@article{KirchStoehr2022,
  author  = {Kirch, Claudia and Stoehr, Christina},
  title   = {Sequential Change Point Tests Based on {U}-Statistics},
  journal = {Scandinavian Journal of Statistics},
  year    = {2022},
  volume  = {49},
  number  = {3},
  pages   = {1184--1214},
  doi     = {10.1111/sjos.12558}
}

@article{KoltchinskiiGine2000,
  author  = {Koltchinskii, Vladimir and Gin{\'e}, Evarist},
  title   = {Random Matrix Approximation of Spectra of Integral Operators},
  journal = {Bernoulli},
  year    = {2000},
  volume  = {6},
  number  = {1},
  pages   = {113--167},
  doi     = {10.2307/3318636}
}

@book{KorolyukBorovskich1994,
  author    = {Korolyuk, Vladimir S. and Borovskich, Yuri V.},
  title     = {Theory of {U}-Statistics},
  publisher = {Kluwer Academic Publishers},
  address   = {Dordrecht},
  year      = {1994}
}

@book{Lee1990,
  author    = {Lee, Alan J.},
  title     = {{U}-Statistics: Theory and Practice},
  publisher = {Marcel Dekker},
  address   = {New York},
  year      = {1990}
}

@article{MannWhitney1947,
  author  = {Mann, Henry B. and Whitney, Donald R.},
  title   = {On a Test of Whether One of Two Random Variables Is Stochastically Larger Than the Other},
  journal = {The Annals of Mathematical Statistics},
  year    = {1947},
  volume  = {18},
  number  = {1},
  pages   = {50--60},
  doi     = {10.1214/aoms/1177730491}
}

@article{MukhopadhyayVik1985,
  author  = {Mukhopadhyay, Nitis and Vik, Inger},
  title   = {Asymptotic Results for Stopping Times Based on {U}-Statistics},
  journal = {Sequential Analysis},
  year    = {1985},
  volume  = {4},
  number  = {1--2},
  pages   = {83--109},
  doi     = {10.1080/07474948508836073}
}

@article{MukhopadhyayVik1988,
  author  = {Mukhopadhyay, Nitis and Vik, Inger},
  title   = {Convergence Rates for Two-Stage Confidence Intervals Based on {U}-Statistics},
  journal = {Annals of the Institute of Statistical Mathematics},
  year    = {1988},
  volume  = {40},
  number  = {1},
  pages   = {111--117}
}

@misc{Nasari2009,
  author       = {Nasari, Masoud M.},
  title        = {Studentized Processes of {U}-Statistics},
  year         = {2009},
  eprint       = {0906.5101},
  archivePrefix= {arXiv},
  primaryClass = {math.ST}
}

@article{RamdasGrunwaldVovkShafer2023,
  author  = {Ramdas, Aaditya and Gr{\"u}nwald, Peter and Vovk, Vladimir and Shafer, Glenn},
  title   = {Game-Theoretic Statistics and Safe Anytime-Valid Inference},
  journal = {Statistical Science},
  year    = {2023},
  volume  = {38},
  number  = {4},
  pages   = {576--601},
  doi     = {10.1214/23-STS894}
}

@article{RosascoBelkinDeVito2010,
  author  = {Rosasco, Lorenzo and Belkin, Mikhail and De Vito, Ernesto},
  title   = {On Learning with Integral Operators},
  journal = {Journal of Machine Learning Research},
  year    = {2010},
  volume  = {11},
  pages   = {905--934}
}

@article{Sproule1985,
  author  = {Sproule, Raymond N.},
  title   = {Sequential Nonparametric Fixed-Width Confidence Intervals for {U}-Statistics},
  journal = {The Annals of Statistics},
  year    = {1985},
  volume  = {13},
  number  = {1},
  pages   = {228--235},
  doi     = {10.1214/aos/1176346588}
}

@article{SzekelyRizzo2013,
  author  = {Sz{\'e}kely, G{\'a}bor J. and Rizzo, Maria L.},
  title   = {Energy Statistics: A Class of Statistics Based on Distances},
  journal = {Journal of Statistical Planning and Inference},
  year    = {2013},
  volume  = {143},
  number  = {8},
  pages   = {1249--1272},
  doi     = {10.1016/j.jspi.2013.03.018}
}

@book{Ville1939,
  author    = {Ville, Jean},
  title     = {{\'{E}}tude Critique de la Notion de Collectif},
  publisher = {Gauthier-Villars},
  address   = {Paris},
  year      = {1939}
}

@article{WaudbySmithArbourSinhaKennedyRamdas2024,
  author  = {Waudby-Smith, Ian and Arbour, David and Sinha, Ritwik and Kennedy, Edward H. and Ramdas, Aaditya},
  title   = {Time-Uniform Central Limit Theory and Asymptotic Confidence Sequences},
  journal = {The Annals of Statistics},
  year    = {2024},
  volume  = {52},
  number  = {6},
  pages   = {2613--2640},
  doi     = {10.1214/24-AOS2408}
}

@article{ArconesGine1993,
  author  = {Arcones, Miguel A. and Gin{\'e}, Evarist},
  title   = {Limit Theorems for {U}-Processes},
  journal = {The Annals of Probability},
  year    = {1993},
  volume  = {21},
  number  = {3},
  pages   = {1494--1542},
  doi     = {10.1214/aop/1176989128}
}

@article{BaringhausFranz2004,
  author  = {Baringhaus, Ludwig and Franz, Carsten},
  title   = {On a New Multivariate Two-Sample Test},
  journal = {Journal of Multivariate Analysis},
  year    = {2004},
  volume  = {88},
  number  = {1},
  pages   = {190--206},
  doi     = {10.1016/S0047-259X(03)00079-4}
}

@incollection{DehlingFriedGarciaWendler2015,
  author    = {Dehling, Herold and Fried, Roland and Garc{\'i}a, Isabel and Wendler, Martin},
  title     = {Change-Point Detection under Dependence Based on Two-Sample {U}-Statistics},
  booktitle = {Asymptotic Laws and Methods in Stochastics},
  series    = {Fields Institute Communications},
  volume    = {76},
  publisher = {Springer},
  address   = {New York},
  year      = {2015},
  pages     = {195--220},
  doi       = {10.1007/978-1-4939-3076-0_12}
}

@article{DehlingVukWendler2022,
  author  = {Dehling, Herold and Vuk, Kata and Wendler, Martin},
  title   = {Change-Point Detection Based on Weighted Two-Sample {U}-Statistics},
  journal = {Electronic Journal of Statistics},
  year    = {2022},
  volume  = {16},
  number  = {1},
  pages   = {862--891},
  doi     = {10.1214/21-EJS1964}
}

@article{SejdinovicSriperumbudurGrettonFukumizu2013,
  author  = {Sejdinovic, Dino and Sriperumbudur, Bharath and Gretton, Arthur and Fukumizu, Kenji},
  title   = {Equivalence of Distance-Based and {RKHS}-Based Statistics in Hypothesis Testing},
  journal = {The Annals of Statistics},
  year    = {2013},
  volume  = {41},
  number  = {5},
  pages   = {2263--2291},
  doi     = {10.1214/13-AOS1140}
}

@article{SzekelyRizzoBakirov2007,
  author  = {Sz{\'e}kely, G{\'a}bor J. and Rizzo, Maria L. and Bakirov, Nail K.},
  title   = {Measuring and Testing Dependence by Correlation of Distances},
  journal = {The Annals of Statistics},
  year    = {2007},
  volume  = {35},
  number  = {6},
  pages   = {2769--2794},
  doi     = {10.1214/009053607000000505}
}

@article{ChoiMarden1998,
  author  = {Choi, Kyungmee and Marden, John I.},
  title   = {A Multivariate Version of Kendall's Tau},
  journal = {Journal of Nonparametric Statistics},
  year    = {1998},
  volume  = {9},
  number  = {3},
  pages   = {261--293},
  doi     = {10.1080/10485259808832746}
}

@inproceedings{ChwialkowskiRamdasSejdinovicGretton2015,
  author    = {Chwialkowski, Kacper and Ramdas, Aaditya and Sejdinovic, Dino and Gretton, Arthur},
  title     = {Fast Two-Sample Testing with Analytic Representations of Probability Measures},
  booktitle = {Advances in Neural Information Processing Systems 28},
  year      = {2015}
}

@article{sgaU2026,
  title={Strong {G}aussian approximation for {U}-statistics in high dimensions and beyond}, 
  author={Li, Weijia and Cai, Leheng and Hu, Qirui},
  journal={arXiv preprint arXiv:2603.10595},
  year={2026},
  eprint={2603.10595},
  archivePrefix={arXiv},
  primaryClass={math.ST}
}

@article{DuerreTylerVogel2016,
  author  = {D{\"u}rre, Alexander and Tyler, David E. and Vogel, Daniel},
  title   = {On the Eigenvalues of the Spatial Sign Covariance Matrix in More Than Two Dimensions},
  journal = {Statistics \& Probability Letters},
  year    = {2016},
  volume  = {111},
  pages   = {80--85},
  doi     = {10.1016/j.spl.2016.01.009}
}

@article{HanLiu2018,
  author  = {Han, Fang and Liu, Han},
  title   = {{ECA}: High-Dimensional Elliptical Component Analysis in Non-Gaussian Distributions},
  journal = {Journal of the American Statistical Association},
  year    = {2018},
  volume  = {113},
  number  = {521},
  pages   = {252--268},
  doi     = {10.1080/01621459.2016.1246366}
}

@inproceedings{ZarembaGrettonBlaschko2013,
  author    = {Zaremba, Wojciech and Gretton, Arthur and Blaschko, Matthew B.},
  title     = {{B}-Test: A Non-Parametric, Low Variance Kernel Two-Sample Test},
  booktitle = {Advances in Neural Information Processing Systems 26},
  year      = {2013},
  pages     = {755--763}
}

}

%%%%%%%%%%%%%%%%%%%%%%%%%%%%%%%%%%%%%%%%%%%%%%%%%%%%%%%%%%%%
\clearpage
\appendix
\setcounter{table}{0}
\setcounter{figure}{0}
\renewcommand{\thetable}{A.%
\arabic{table}}
\renewcommand{\thefigure}{A.%
\arabic{figure}}
\section{Additional numerical results}

\subsection{Results of sensitivity experiments}\label{sec:sensitivity}

\begin{figure*}[h]
    \centering
    \begin{subfigure}[b]{0.48\textwidth}
        \centering
        \includegraphics[width=\textwidth]{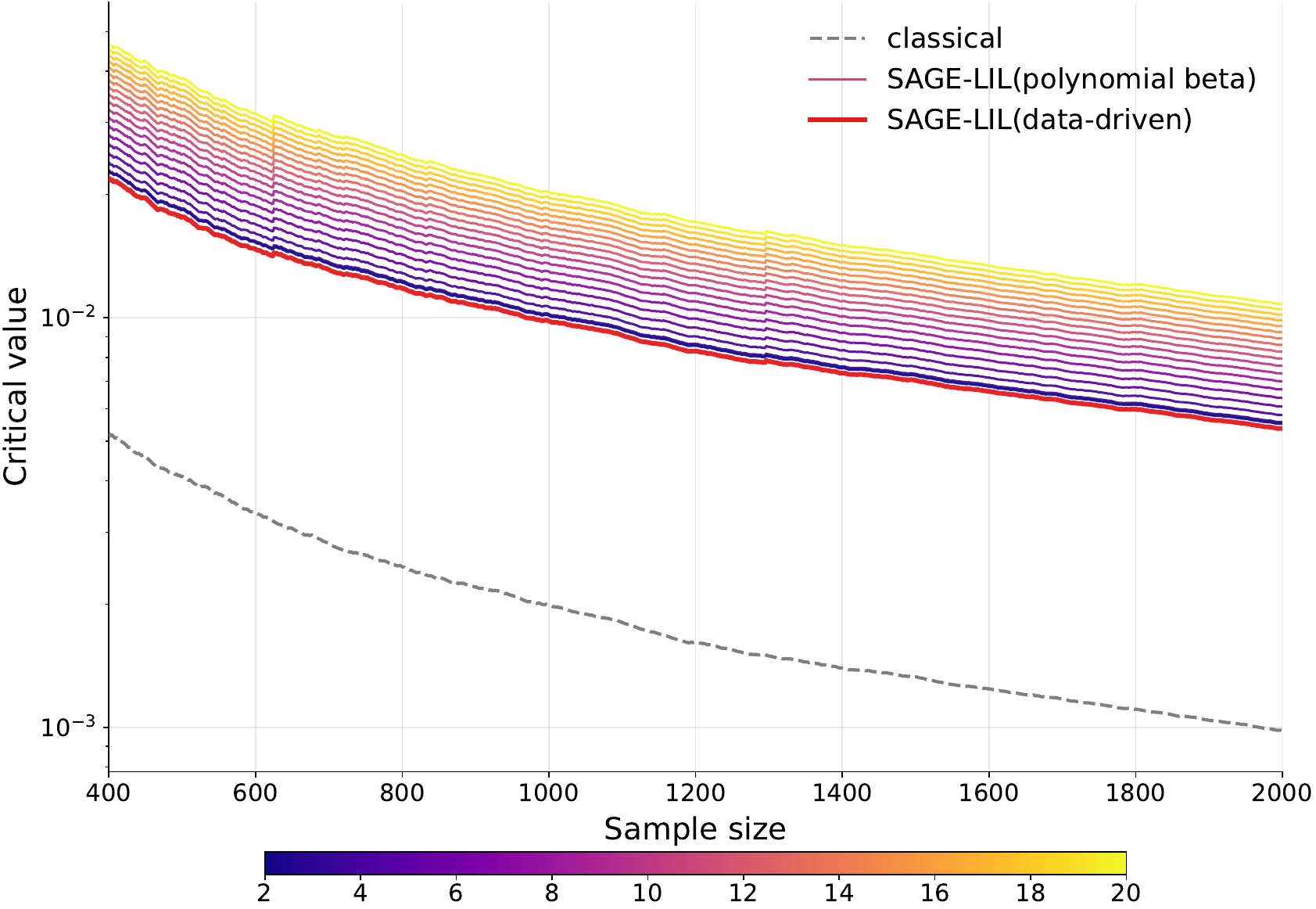} 
        \caption{Polynomial decay ($\beta_\ell \propto \ell^{-b}$)}
        \label{fig:deg_sens_poly}
    \end{subfigure}
    \hfill
    \begin{subfigure}[b]{0.48\textwidth}
        \centering
        \includegraphics[width=\textwidth]{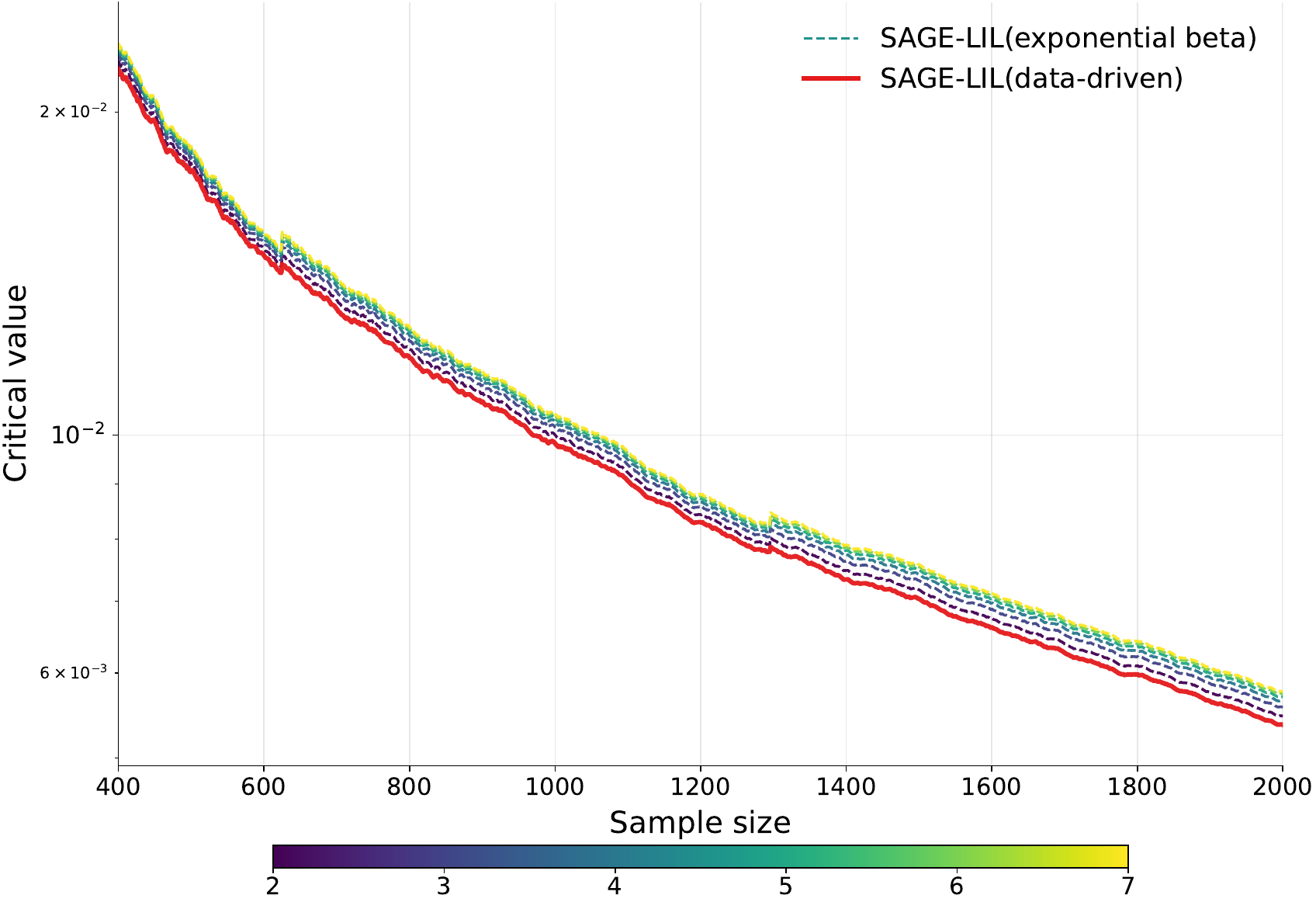} 
        \caption{Exponential decay ($\beta_\ell \propto e^{-c(\ell-1)}$)}
        \label{fig:deg_sens_exp}
    \end{subfigure}
    \caption{Sensitivity of the width of SAGE-LIL to the spectral weight allocation $\beta_\ell$. The colorbars indicate the decay-rate parameters: the polynomial degree \(b \in [2,20]\) in panel (a), and the exponential rate \(c \in [2,7]\) in panel (b).}
    \label{fig:deg_sensitivity}
\end{figure*}

Figure \ref{fig:deg_sensitivity} reports the sensitivity analysis for the weight allocation, as discussed in Section \ref{sec:simulations}. The results show that using predefined deterministic weights, such as polynomial allocations $\beta_\ell \propto \ell^{-b}$ with a large value of $b$, can lead to overly conservative boundaries. To improve efficiency, the data-driven  allocation suggests choosing weights proportional to the positive part of the spectrum, namely
\(
\beta_\ell^\star \propto \max\{\lambda_\ell,0\}
\). 
Since the population spectrum is unknown, in our simulations we adopt the dynamically updated plug-in weights
\(
\widehat{\beta}_{\ell,n} \propto \max\{\widehat{\lambda}_{\ell,n},0\}
\).

We note a minor but important technical distinction between the theoretical SAGE boundary in Theorem \ref{thm:sequential_test} and its practical implementation in the simulations. In the theorem, the weights are required to be a fixed deterministic nonnegative array satisfying $\sum_{\ell\ge1}\beta_{\ell,n}=1$, which ensures the time-uniform union bound
\[
\mathbb{P}\Big(\bigcup_{\ell\ge1} E_\ell^c\Big)
\le
\sum_{\ell\ge1}\alpha\beta_{\ell,n}
=
\alpha.
\]
By contrast, the simulation procedure uses the data-driven weights $\widehat{\beta}_{\ell,n}$, which are random. Thus, this implementation should be viewed as an empirical plug-in heuristic. As shown in Figure \ref{fig:deg_sensitivity}, this heuristic empirically maintains nominal coverage and delivers robust practical performance. Providing a rigorous time-uniform validity theory for SAGE boundaries with dynamically updated random weights is an interesting direction for future work.

\clearpage
\subsection{AsympCSs for GMD under Laplace and $t_{10}$ distributions}\label{app:gmd}
\begin{figure}[h] 
    \centering
    \includegraphics[width=\textwidth]{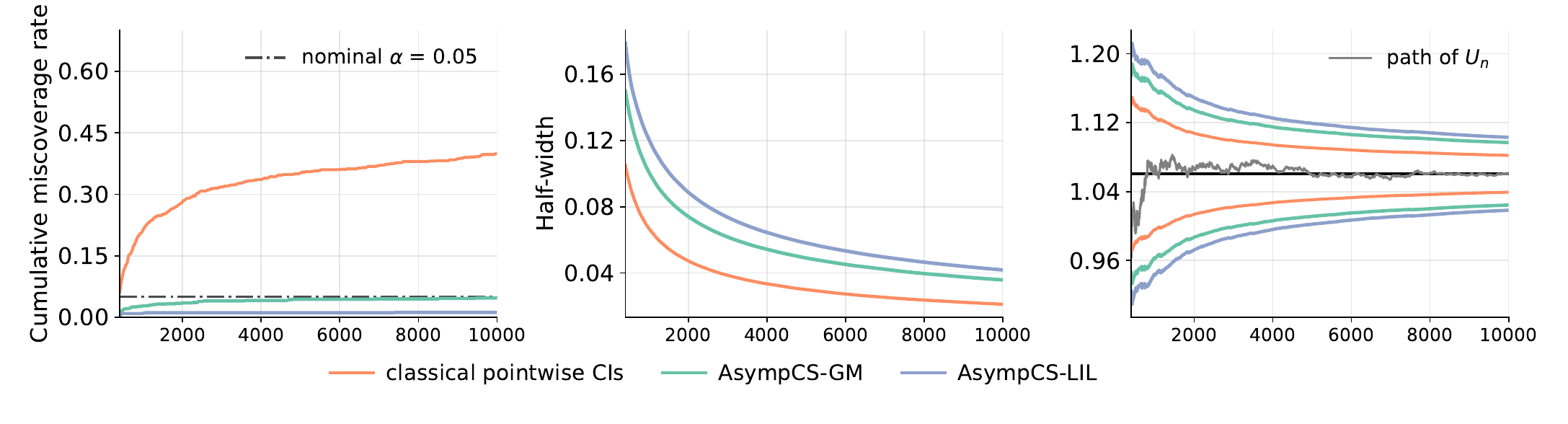} 
    \caption{Left: cumulative miscoverage rates of the proposed AsympCSs and the classical pointwise CIs for GMD under the Laplace distribution. Middle: averaged half-widths of the proposed AsympCSs and the classical pointwise CIs. Right: a single sample path of the statistics $U_n$ for GMD alongside the three boundaries, where the black horizontal line indicates the true parameter. The horizontal axis in all three panels represents the sample size.}
    \label{fig:la-gmd}
\end{figure}

\begin{figure}[h] 
    \centering
    \includegraphics[width=\textwidth]{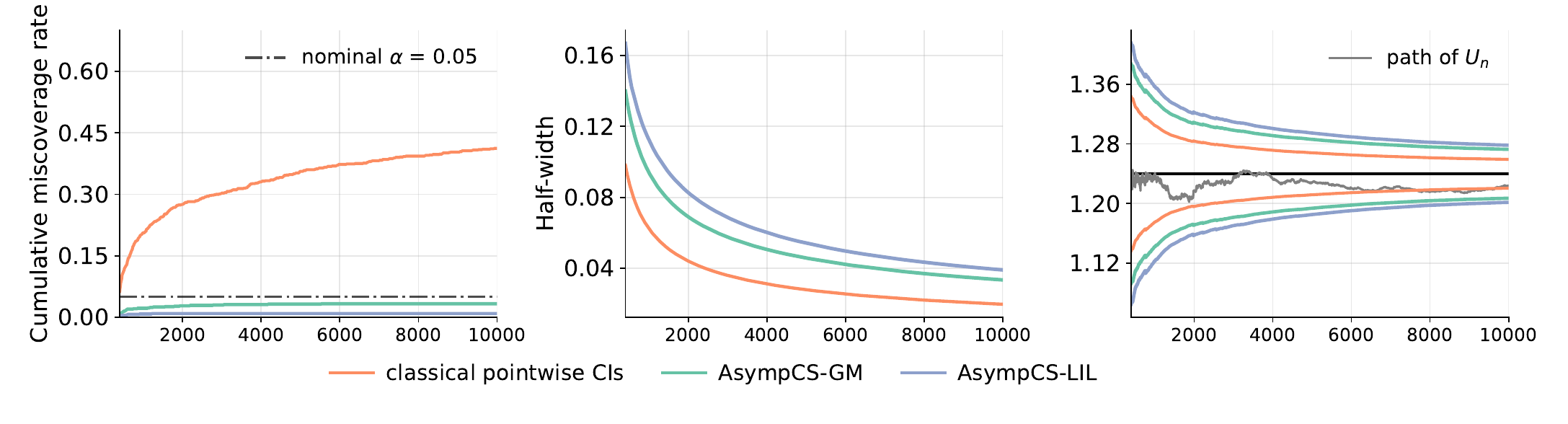} 
    \caption{Left: cumulative miscoverage rates of the proposed AsympCSs and the classical pointwise CIs for GMD under the $t_{10}$ distribution. Middle: averaged half-widths of the proposed AsympCSs and the classical pointwise CIs. Right: a single sample path of the statistics $U_n$ for GMD alongside the three boundaries, where the black horizontal line indicates the true parameter. The horizontal axis in all three panels represents the sample size.}
    \label{fig:t-gmd}
\end{figure}

Figures \ref{fig:la-gmd}-\ref{fig:t-gmd} present the performance of the proposed AsympCSs for GMD under Laplace and $t_{10}$ distributions. {In our implementation, the data dimension is set to $d=1$ for the GMD, sample variance, and MMD tests, and $d=2$ for the Spatial Kendall's tau.}

\clearpage
  
\subsection{Two-sample sequential testing based on MMD  under  Laplace and $t_{10}$ null distributions}\label{app:mmd}
\begin{figure}[h] 
    \centering
    \includegraphics[width=\textwidth]{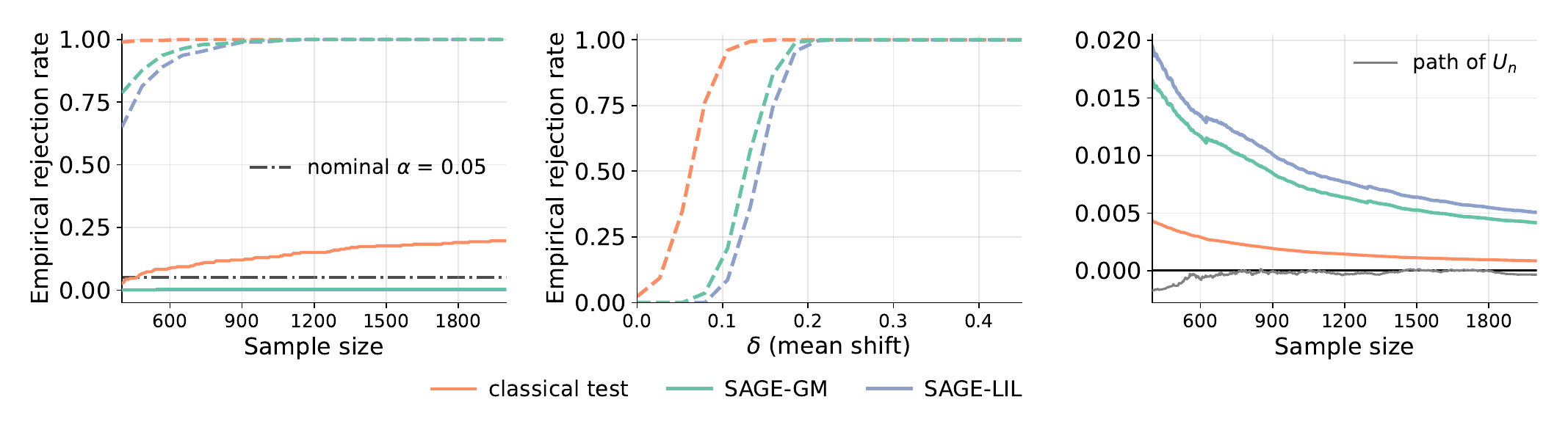} 
    \caption{Left: the size and power comparison of the proposed sequential testing procedure using the SAGE boundaries and the classical test procedure for sequential two-sample test with MMD kernel statistics under Laplace distribution. Dashed lines represent the power under $H_1$ with  $\delta = 0.3$, while solid lines represent the size under $H_0$ with $\delta = 0$. Middle: the empirical power   over mean shifts $\delta \in [0, 0.45]$ at sample size $n = 2000$. Right: a single sample path of the MMD test statistic $U_n$ alongside the three boundaries under $H_0$.}
    \label{fig:la-mmd}
\end{figure}

\begin{figure}[h] 
    \centering
    \includegraphics[width=\textwidth]{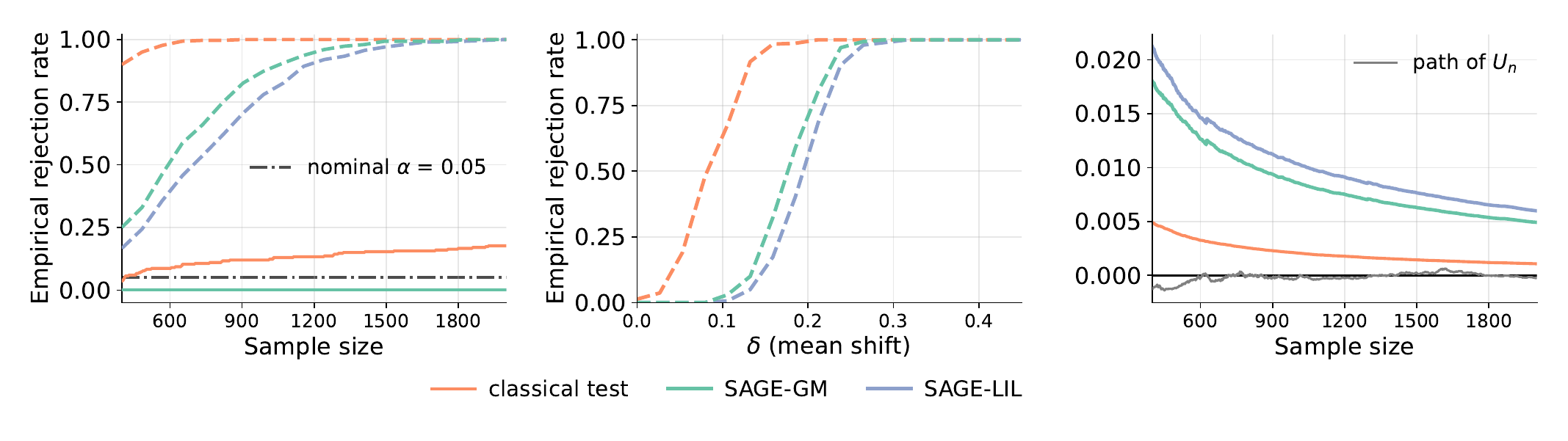} 
    \caption{Left: the size and power comparison of the proposed sequential testing procedure using the SAGE boundaries and the classical test procedure for sequential two-sample test with MMD kernel statistics under $t_{10}$ distribution. Dashed lines represent the power under $H_1$ with  $\delta = 0.3$, while solid lines represent the size under $H_0$ with $\delta = 0$. Middle: the empirical power   over mean shifts $\delta \in [0, 0.45]$ at sample size $n = 2000$. Right: a single sample path of the MMD test statistic $U_n$ alongside the three boundaries under $H_0$.}
    \label{fig:t-mmd}
\end{figure}

Figures \ref{fig:la-mmd}-\ref{fig:t-mmd} present the performance of the sequential testing procedures using
the proposed SAGE boundaries for sequential two-sample test with MMD
kernel statistics under Laplace and $t_{10}$ distributions. 

\clearpage
\subsection{Additional experiments}\label{app:addition-fig}
\begin{figure}[h] 
    \centering
    \includegraphics[width=\textwidth]{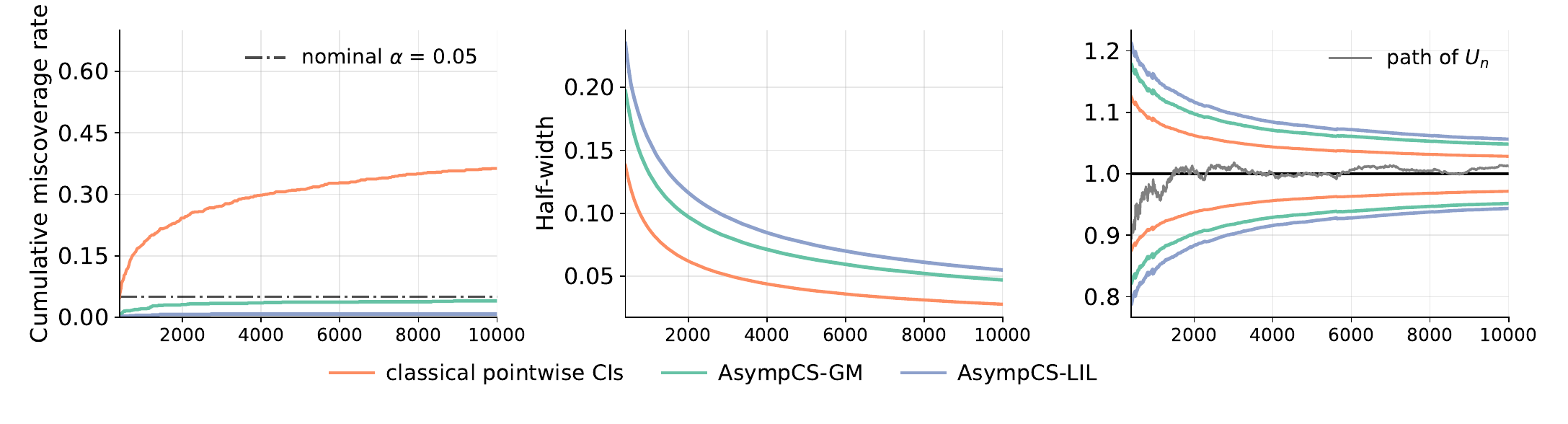} 
    \caption{Left: cumulative miscoverage rates of the proposed AsympCSs and the classical pointwise CIs for sample variance under the Gaussian distribution. Middle: averaged half-widths of the proposed AsympCSs and the classical pointwise CIs. Right: a single sample path of the statistics $U_n$ for sample variance alongside the three boundaries, where the black horizontal line indicates the true parameter. The horizontal axis in all three panels represents the sample size.}    
    \label{fig:gauss-sv}
\end{figure}

\begin{figure}[h] 
    \centering
    \includegraphics[width=\textwidth]{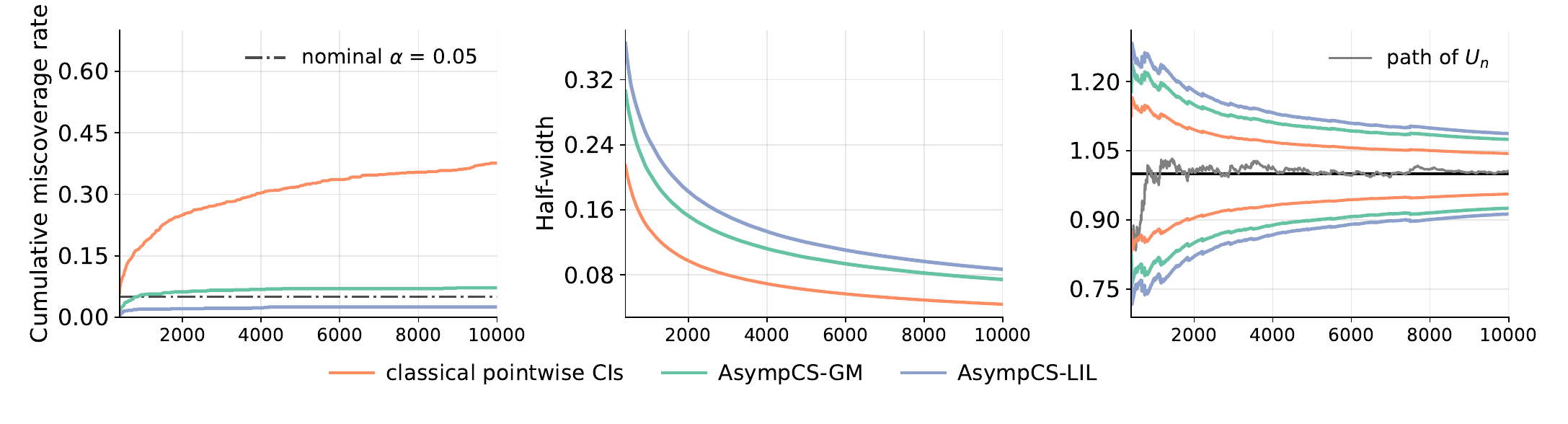} 
    \caption{Left: cumulative miscoverage rates of the proposed AsympCSs and the classical pointwise CIs for sample variance under the Laplace distribution. Middle: averaged half-widths of the proposed AsympCSs and the classical pointwise CIs. Right: a single sample path of the statistics $U_n$ for sample variance alongside the three boundaries, where the black horizontal line indicates the true parameter. The horizontal axis in all three panels represents the sample size.}
    \label{fig:la-sv}
\end{figure}

\begin{figure}[h] 
    \centering
    \includegraphics[width=\textwidth]{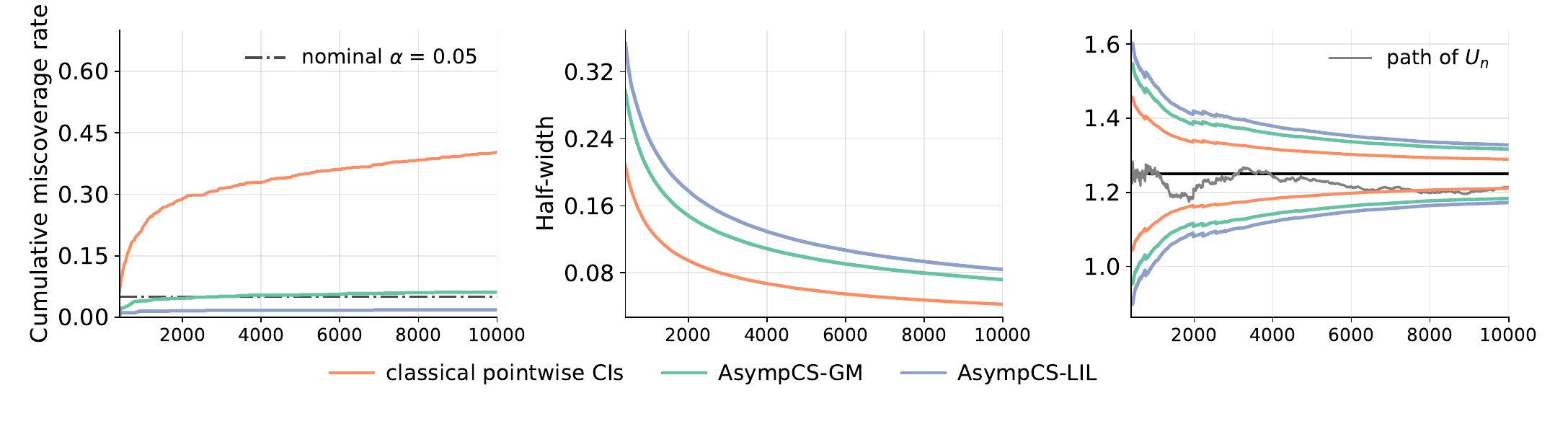} 
    \caption{Left: cumulative miscoverage rates of the proposed AsympCSs and the classical pointwise CIs for sample variance under the $t_{10}$ distribution. Middle: averaged half-widths of the proposed AsympCSs and the classical pointwise CIs. Right: a single sample path of the statistics $U_n$ for sample variance alongside the three boundaries, where the black horizontal line indicates the true parameter. The horizontal axis in all three panels represents the sample size.}
    \label{fig:t-sv}
\end{figure}

Figures \ref{fig:gauss-sv}-\ref{fig:t-sv} present the performance of the proposed AsympCSs for sample variance under the same three distributions as in Section~\ref{sec:simulations}.

We  notice that in Figures \ref{fig:la-sv}-\ref{fig:t-sv}, the cumulative miscoverage rate of AsympCS-GM   slightly exceeds the nominal level. Theoretically,   the Gaussian approximation errors for \(U_n\) and the consistency of the variance estimator $\widehat\sigma_n^2$ depend on higher-order moments of the kernel. For the sample variance, such moment dependence is more pronounced, so heavier-tailed distributions require larger sample sizes for the asymptotic approximation to become accurate. 
Increasing the cold-start parameter \(m\) mitigates this finite-sample under-coverage by delaying monitoring until sufficient data have accumulated. 

\begin{figure}[h] 
    \centering
    \includegraphics[width=\textwidth]{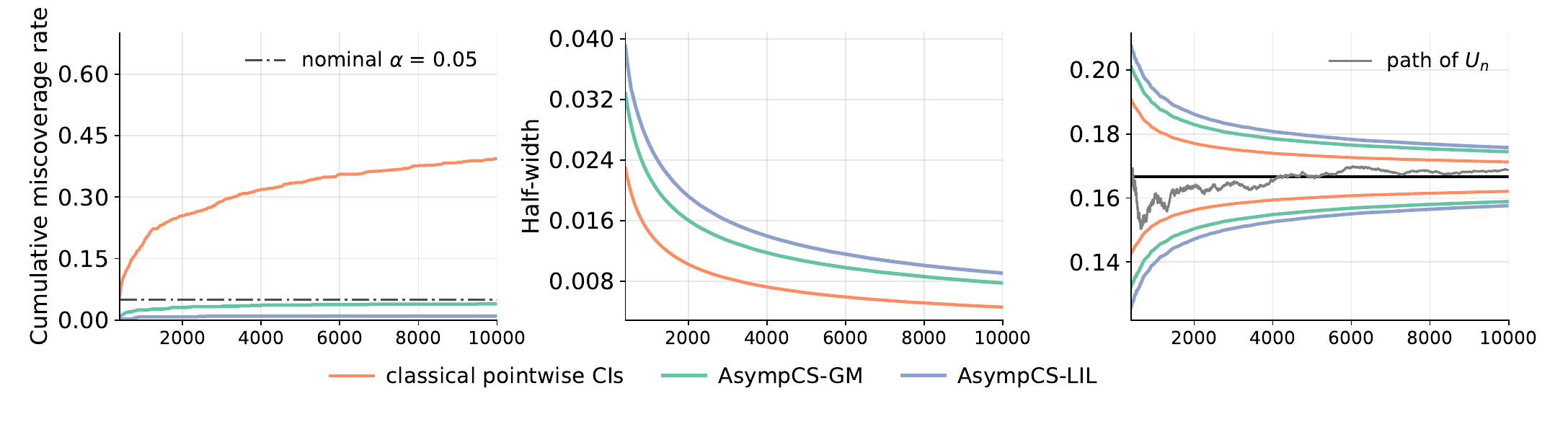} 
    \caption{Left: cumulative miscoverage rates of the proposed AsympCSs and the classical pointwise CIs for Spatial Kendall's tau under the Gaussian distribution. Middle: averaged half-widths of the proposed AsympCSs and the classical pointwise CIs. Right: a single sample path of the statistics $U_n$ for Spatial Kendall's tau alongside the three boundaries, where the black horizontal line indicates the true parameter. The horizontal axis in all three panels represents the sample size.}
    \label{fig:gauss-sk}
\end{figure}

\begin{figure}[h] 
    \centering
    \includegraphics[width=\textwidth]{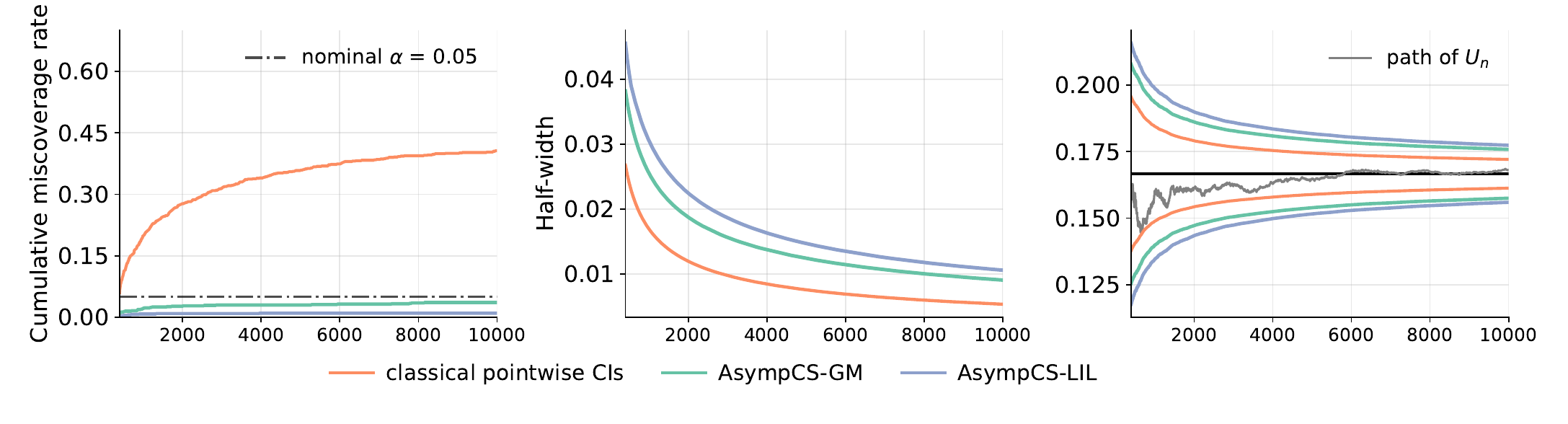} 
    \caption{Left: cumulative miscoverage rates of the proposed AsympCSs and the classical pointwise CIs for Spatial Kendall's tau under the Laplace distribution. Middle: averaged half-widths of the proposed AsympCSs and the classical pointwise CIs. Right: a single sample path of the statistics $U_n$ for Spatial Kendall's tau alongside the three boundaries, where the black horizontal line indicates the true parameter. The horizontal axis in all three panels represents the sample size.}
    \label{fig:la-sk}
\end{figure}

\begin{figure}[h] 
    \centering
    \includegraphics[width=\textwidth]{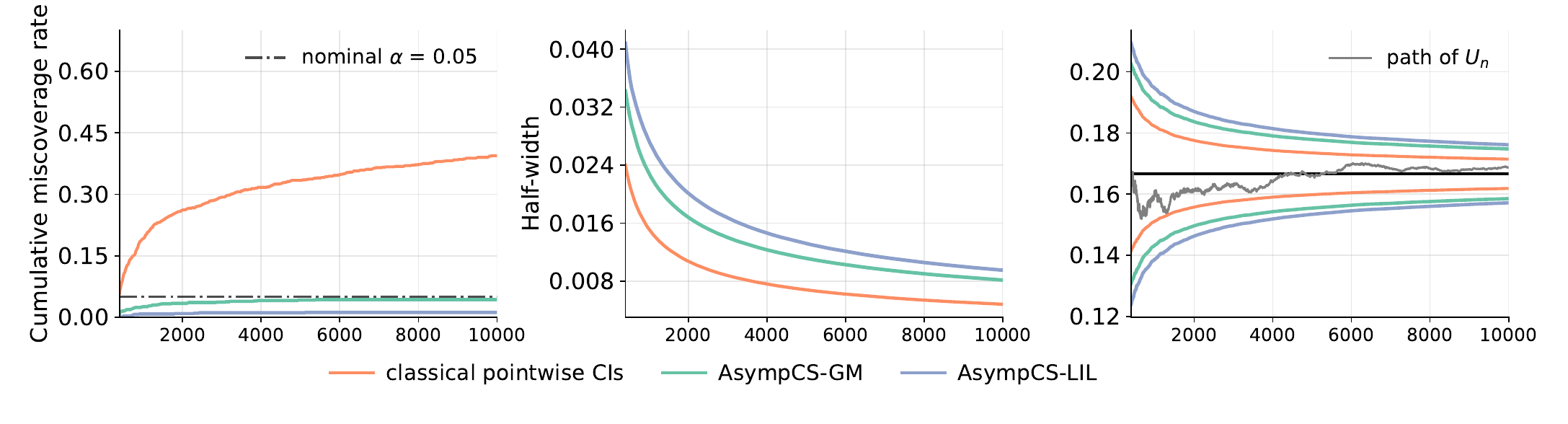} 
    \caption{Left: cumulative miscoverage rates of the proposed AsympCSs and the classical pointwise CIs for Spatial Kendall's tau under the $t_{10}$ distribution. Middle: averaged half-widths of the proposed AsympCSs and the classical pointwise CIs. Right: a single sample path of the statistics $U_n$ for Spatial Kendall's tau alongside the three boundaries, where the black horizontal line indicates the true parameter. The horizontal axis in all three panels represents the sample size.}
    \label{fig:t-sk}
\end{figure}

Figures \ref{fig:gauss-sk}-\ref{fig:t-sk} present the performance of the proposed AsympCSs estimation for Spatial Kendall’s tau under the three distributions. 
Here, we detail the data generation process for the Spatial Kendall's tau estimation. Unlike the other three kernels we adopted, which  utilize one-dimensional data drawn from the three specified   distributions, Spatial Kendall's tau requires bivariate observations.  To accommodate this, we generate two-dimensional data vectors $\boldsymbol{X}_i = (X_{i1}, X_{i2})^\top \in \mathbb{R}^2$ with a  correlation structure. Specifically, the dependence is governed by a shape matrix $\mathbf{\Sigma}$, defined as:
\begin{align*}
    \mathbf{\Sigma} = \begin{pmatrix} 1 & \rho \\ \rho & 1 \end{pmatrix}.
\end{align*}
Throughout our simulations, the shape correlation is set to $\rho = 0.6$. For the two-dimensional elliptical distributions, $\bm{X}_i$ has the scale mixture of normals representation:
\begin{equation}
    \bm{X}_i = \sqrt{W_i} \mathbf{A} \bm{Z}_i,
    \label{eq:ellip-data}
\end{equation}

where $\boldsymbol{Z}_i \sim \mathcal{N}(\mathbf{0}, \mathbf{I}_2)$ is a standard two-dimensional Gaussian vector, and $\mathbf{A}$ is the Cholesky factor satisfying $\mathbf{A}\mathbf{A}^\top = \boldsymbol{\Sigma}$. The   multiplier  $W_i > 0$ is a scalar random variable, independent of $\boldsymbol{Z}_i$, that dictates the tail heaviness of the resulting bivariate distribution. 

In our simulations, we use three different $W_i$:
\begin{itemize}[leftmargin=*]
    \item \textbf{Gaussian}: $W_i = 1$ a.s., yielding the  bivariate normal distribution with correlation $\rho$.
    \item \textbf{Student's $t_{10}$}: $W_i = \nu / V_i$, where $V_i \sim \chi^2_\nu$ with degrees of freedom $\nu=10$.
    \item \textbf{Laplace}: $W_i \sim \mathrm{Exp}(1)$, which generates the bivariate symmetric Laplace distribution.
\end{itemize}
Based on this representation, we first sample the standard two-dimensional Gaussian data $\bm{Z}_i$ and obtain $\bm{X}_i$ from \eqref{eq:ellip-data}, then compute the corresponding  U-statistics. By mathematical derivation,  
\begin{align*}
    \theta   =\mathbb{E}\left\{\frac{({X}_{11}-{X}_{21})({X}_{12}-{X}_{22})}{\|\bm{X}_{1}-\bm{X}_{2}\|_2^2} \mathbf{1}\{\bm X_1 \neq \bm X_2\} \right\}= \left(1-\sqrt{1-\rho^2}\right)/(2\rho).
\end{align*}

\clearpage

\subsection{Discussion of  computational complexity}\label{app:comp-complexity}

Algorithm \ref{alg:degenerate} requires an eigen-decomposition of an \(n\times n\) Gram matrix at each monitoring time, leading to a per-step computational cost of \(\mathcal O(n^3)\). To reduce this cost, one may approximate the eigenvalues using an \(N_n\times N_n\) sub-Gram matrix formed from the first
\(
N_n=\lceil n^\varpi\rceil
\)
observations, where \(\varpi\in(0,1)\). For example, taking \(\varpi=2/3\) reduces the eigen-decomposition cost to \(\mathcal O(N_n^3)=\mathcal O(n^2)\). 
Under the empirical  estimation of eigenvalues discussed in Remark \ref{remark1}, this subsampling strategy still preserves the asymptotic anytime-validity of the resulting plug-in SAGE boundary, provided that the truncation number $L_n\asymp n^{a}$ with \(0<a<1/3\).

\section{Further discussion}
\subsection{Discussion on  the lower boundaries}\label{sec:lower}

Recall that $\Lambda = \sum_{\ell\geq1 } \lambda_\ell $. Define $\Lambda^- = \sum_{\ell:\lambda_\ell<0 } \lambda_\ell$, $\Lambda_\beta^- = \sum_{\ell:\lambda_\ell<0 } \lambda_\ell \log(1/\beta_l)$ and $\Lambda_{\beta,g}^-=\sum_{\ell:\lambda_\ell<0} \lambda_\ell  \left\{g^{-1}(\alpha\beta_\ell)\right\}^2$,  where \(g(a)=2\{1-\Phi(a)+a\phi(a)\}\). We now introduce the  lower boundary $\mathcal{L}_{\alpha,m}(\cdot)$.  
Suppose that $\sum_{\ell\geq 1}|\lambda_\ell|<\infty$, $\Lambda_\beta^-<\infty$  and \(\{ W_\ell(\cdot)\}_{\ell\ge 1}\) are independent standard Brownian motions.  
For any \(\alpha\in(0,1)\), \(m\ge 1\), 
\[
\mathbb{P}\left(
\forall n\ge m:\ 
\frac{1}{n^2}\sum_{\ell\ge 1}\lambda_\ell\bigl\{W_\ell^2(n)-n\bigr\}
\geq
\mathcal{L}_{\alpha,m}(n)
\right)
\ge 1-\alpha,
\]
where \(\mathcal{L}_{\alpha,m}(n)\) can be chosen in either of the following two ways:
    \begin{align*}%\label{LIL2}
        &\mathcal{L}_{\alpha,m}^{\text{LIL}}(n)=\frac{\left(\eta^{1/4}+\eta^{-1/4}\right)^2}{ 2n }\left[  \left\{ {s\log\log\left(\max\left\{ \frac{\eta n}{m},e\right\}\right)+\log \frac{\zeta(s)}{\alpha (\log\eta)^s}  } \right\}\Lambda^- + \Lambda_\beta^- \right]-\frac{\Lambda}{n},
    \end{align*}  
\begin{align*}%\label{GM2}
    &\mathcal{L}_{\alpha,m}^{\text{GM}}(n)=\frac{1}{n}\left[  \Lambda^-\log\left(\frac{n}{m}\right)+\Lambda_{\beta,g}^--\Lambda\right].
\end{align*}
If the kernel is positive semi-definite, then all eigenvalues are nonnegative, so the SAGE lower boundary degenerates to the deterministic term \(\mathcal{L}_{\alpha,m}^{\text{LIL}}(n)=\mathcal{L}_{\alpha,m}^{\text{GM}}(n)=-\Lambda/n\).

\subsection{Discussion on Remark \ref{remark1}}\label{sec:details}
Let \(\mathcal K:L^2(\mathbb P)\to L^2(\mathbb P)\) be the integral operator induced by the centered kernel
\[
K(x,y):=h(x,y)-\theta,
\]
that is,
\[
(\mathcal K f)(x)
=
\int K(x,y)f(y)\,d\mathbb P(y)
=
\E\{[h(x,X)-\theta]f(X)\}.
\]
Let \((\lambda_\ell)_{\ell\ge1}\) denote the eigenvalues of \(\mathcal K\).

Given the sample \(X_1,\ldots,X_n\), define the empirical operator \(\widetilde{\mathcal K}_n\) by
\[
(\widetilde{\mathcal K}_n f)(x)
=
\frac1n\sum_{j=1}^n
\{h(x,X_j)-\theta\}f(X_j).
\]
Equivalently, when restricted to the sample points, \(\widetilde{\mathcal K}_n\) is represented by the matrix
\(n^{-1}\widetilde{\mathbf K}_n\), where
\[
\widetilde {  K}_n(i,j):=h(X_i,X_j)-\theta,
\qquad 1\le i,j\le n.
\]
Let \((\widetilde\lambda_{\ell,n})_{\ell\ge1}\) be the eigenvalues of
\(n^{-1}\widetilde{\mathbf K}_n\).

Similarly, consider the centered Gram matrix
\[
\widehat K_n(i,j)
:=
h(X_i,X_j)-U_n,
\qquad 1\le i,j\le n,
\]
and let \((\widehat\lambda_{\ell,n})_{\ell\ge1}\) be the eigenvalues of
\(n^{-1}\widehat{\mathbf K}_n\). Then
\[
 n^{-1}\widehat {\mathbf K}_n
=
 n^{-1}\widetilde {\mathbf K}_n 
-
n^{-1}(U_n-\theta)\mathbf J_n,
\]
where \(\mathbf J_n\) denotes the \(n\times n\) all-ones matrix.
Since \(n^{-1}\mathbf J_n\) has eigenvalues \(1,0,\ldots,0\),
\[
\left\|n^{-1}(U_n-\theta)\mathbf J_n\right\|_{\mathrm{op}}
=
|U_n-\theta|
=
\O_{a.s.}\!\left(\frac{\log\log n}{n}\right).
\]
By Weyl's inequality,
\[
\sup_{\ell\ge1}
\left|
\widehat\lambda_{\ell,n}
-
\widetilde\lambda_{\ell,n}
\right|
\le
|U_n-\theta|
=
\O_{a.s.}\!\left(\frac{\log\log n}{n}\right).
\]
Moreover, by the results in \cite{KoltchinskiiGine2000,RosascoBelkinDeVito2010} under standard regularity conditions on the kernel function $h(\cdot,\cdot)$, 
\[
\left\|
\widetilde{\mathcal K}_n-\mathcal K
\right\|_{\mathrm{op}}
=
\O_{a.s.}\!\left(\sqrt{\frac{\log n}{n}}\right).
\]
Again by Weyl's inequality,
\[
\sup_{\ell\ge1}
\left|
\widetilde\lambda_{\ell,n}-\lambda_\ell
\right|
\le
\left\|
\widetilde{\mathcal K}_n-\mathcal K
\right\|_{\mathrm{op}}
=
\O_{a.s.}\!\left(\sqrt{\frac{\log n}{n}}\right).
\]
Consequently,
\[
\sup_{\ell\ge1}
\left|
\widehat\lambda_{\ell,n}-\lambda_\ell
\right|
=
\O_{a.s.}\!\left(\sqrt{\frac{\log n}{n}}\right).
\]

Let \(x_+=\max\{x,0\}\) and define $\Lambda^+=\sum_{\ell\ge1}(\lambda_\ell)_+$ and $\Lambda_\beta^+=\sum_{\ell\ge1}\log(1/\beta_\ell)(\lambda_\ell)_+$, 
and their truncated empirical counterparts
$\widehat\Lambda^+
=
\sum_{\ell=1}^{L_n}(\widehat\lambda_\ell)_+$, 
 $
\widehat\Lambda_\beta^+
=
\sum_{\ell=1}^{L_n}\log(1/\beta_\ell)(\widehat\lambda_\ell)_+$.

A sufficient condition for
\[
 \widehat\Lambda^+-\Lambda^+ 
=
\o_{a.s.}(n^{-\iota}),
\quad
 \widehat\Lambda_\beta^+-\Lambda_\beta^+ 
=
\o_{a.s.}(n^{-\iota}), \quad
 \widehat\Lambda_{\beta,g}^+-\sum_{\ell:\lambda_\ell>0} \lambda_\ell  \left\{g^{-1}(\alpha\beta_\ell)\right\}^2 
=
\o_{a.s.}(n^{-\iota}),
\]
is 
\begin{align}\label{sufficient}
    \sum_{\ell=1}^{L_n}
\log(1/\beta_\ell)|\widehat\lambda_\ell-\lambda_\ell|
=
\o_{a.s.}(n^{-\iota}),
\quad
\sum_{\ell>L_n}\log(1/\beta_\ell)(\lambda_\ell)_+
=
\o(n^{-\iota}).
\end{align}
In fact, since \(x\mapsto x_+\) is \(1\)-Lipschitz,
\[
|\widehat\Lambda^+-\Lambda^+|
\le
\sum_{\ell=1}^{L_n}|\widehat\lambda_\ell-\lambda_\ell|
+
\sum_{\ell>L_n}(\lambda_\ell)_+,
\]
and similarly,
\[
|\widehat\Lambda_\beta^+-\Lambda_\beta^+|
\le
\sum_{\ell=1}^{L_n}\log(1/\beta_\ell)|\widehat\lambda_\ell-\lambda_\ell|
+
\sum_{\ell>L_n}\log(1/\beta_\ell)(\lambda_\ell)_+.
\]
Besides, by the tail behavior established in the proof of Theorem \ref{thm:gaussian_chaos}, \begin{align*}
    &\left| \sum_{\ell:\lambda_\ell>0 } \lambda_\ell \{g^{-1}(\alpha\beta_l)\}^2- \sum_{\ell:\widehat\lambda_\ell>0 }\widehat\lambda_\ell \{g^{-1}(\alpha\beta_l)\}^2\right|
 \\&\leq \sum_{\ell>L_n,\lambda_\ell>0} \lambda_\ell \{g^{-1}(\alpha\beta_l)\}^2+
\sum_{\ell=1}^{L_n}|\widehat\lambda_{\ell}-\lambda_\ell|\{g^{-1}(\alpha\beta_l)\}^2 \\&\lesssim  \sum_{\ell>L_n,\lambda_\ell>0} \lambda_\ell\log(1/\beta_\ell)+
\sum_{\ell=1}^{L_n}|\widehat\lambda_{\ell}-\lambda_\ell|\log(1/\beta_\ell). 
\end{align*}
% When $\beta_\ell\asymp \ell^{-b}$ for some $b>1$, $L_n\asymp n^{a}$ for some $0<a<1/2$, we could verify that \eqref{sufficient} holds.
{\color{black}For example, if $\beta_\ell\asymp \ell^{-b}$ with $b>1$ and $L_n\asymp n^{a}$, the leading term is controlled by
\[
\sum_{\ell=1}^{L_n}\log(1/\beta_\ell)|\widehat\lambda_\ell-\lambda_\ell|
=O_{a.s.}\!\left(n^a(\log n)^{3/2}n^{-1/2}\right),
\]
under the displayed operator-norm rate. Hence, $0<a<1/2$ ensures the condition in \eqref{sufficient} holds for some small $0<\iota<1/2-a$.}

\end{document}